\documentclass[12pt]{article}
\usepackage{amsmath}
\usepackage{amssymb}
\usepackage{latexsym}
\usepackage{amsthm}
\usepackage{mathrsfs}
\usepackage{fancyhdr}
\usepackage{amsfonts}
\usepackage{amssymb}%Use Miscell. symbols
\usepackage{graphicx}

\newtheorem{thm}{Theorem}[section]
\newtheorem{lem}[thm]{Lemma}

\newtheorem{pro}[thm]{Proposition}
\newtheorem{defi}[thm]{Definition}

\newtheorem{obs}{Observation}

\usepackage{latexsym , bm}

\title{Extremal fullerene graphs with the maximum Clar number
\thanks{This paper is supported by NSFC grant 10831001.}}
\author{ Dong Ye \ and \ Heping Zhang
\\\small{School of Mathematics and Statistics, Lanzhou University, Lanzhou, Gansu 730000, P.R. China}
\\\small{E-mails: dye@lzu.edu.cn, zhanghp@lzu.edu.cn}}
\date{}

\begin{document}

\maketitle

\begin{abstract}
A fullerene graph is a cubic 3-connected plane graph with (exactly
12) pentagonal faces and hexagonal faces. Let $F_n$ be a fullerene
graph with $n$ vertices. A set $\mathcal H$ of mutually disjoint
hexagons of $F_n$ is a sextet pattern if $F_n$ has a perfect
matching which alternates on and off each hexagon in $\mathcal H$.
The maximum cardinality of  sextet patterns  of $F_n$ is the Clar
number of $F_n$. It was shown that the Clar number is no more than
$\lfloor\frac {n-12} 6\rfloor$. Many fullerenes with experimental
evidence attain the upper bound, for instance, $\text{C}_{60}$ and
$\text{C}_{70}$. In this paper, we characterize extremal fullerene
graphs whose Clar numbers equal $\frac{n-12} 6$. By the
characterization, we show that there are precisely 18 fullerene
graphs with 60 vertices, including $\text{C}_{60}$, achieving the
maximum Clar number 8 and we construct all these extremal fullerene
graphs.

\vspace{0.3cm}

\noindent {\em Keywords:} Fullerene graph; Clar number; Perfect
matching; Sextet pattern; $\text{C}_{60}$

\noindent {\em AMS 2000 subject classification:} 05C70, 05C90
\end{abstract}

%%introduction--------------------------------------------------
\section{Introduction}

A {\em fullerene graph} is a cubic 3-connected plane graph which has
exactly 12 pentagonal faces and other hexagonal faces. Fullerene
graphs correspond to the fullerene molecule frames in chemistry. Let
$F_n$ be a fullerene graph with $n$ vertices. It is well known that
$F_n$ exists for any even $n\ge 20$ except $n=22$ \cite{BD,FM}. For
small $n$, a constructive enumeration of fullerene isomers with $n$
vertices was given \cite{BD}. For example, there are 1812 distinct
fullerene graphs with 60 vertices including the famous
$\text{C}_{60}$ synthesized in 1985 by Kroto et al. \cite{KHOCS}.

\begin{figure}[!hbtp]\refstepcounter{figure}\label{fig1-1}
\begin{center}
\includegraphics{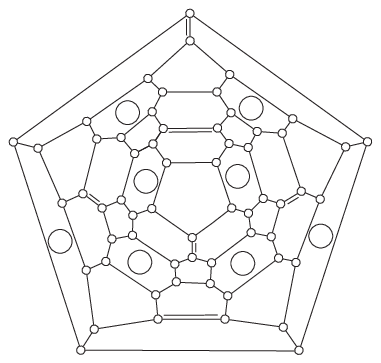}\\
{Figure \ref{fig1-1}: A Clar formula of $\text{C}_{60}$.}
\end{center}
\end{figure}

Let $F$ be a fullerene graph. A {\em perfect matching} (Kekul\'e
structure in chemistry) of $F$ is a set $M$ of independent edges
such that every vertex of $F$ is incident with an edge in $M$. A
cycle of $F$ is {\em $M$-alternating} (or {\em conjugated}) if its
edges appear alternately in and off $M$. A set $\mathcal H$ of
mutually disjoint hexagons is called a {\em sextet pattern} if $F$
has a perfect matching $M$ such that every hexagon in $\mathcal H$
is $M$-alternating. So if $\mathcal H$ is a sextet pattern of $F$,
then $F-\mathcal H$ has a perfect matching where $F-\mathcal H$ is
the subgraph arising from $F$ by deleting all vertices and edges
incident with hexagons in $\mathcal H$. A maximum sextet pattern is
also called a {\em Clar formula}. The cardinality of a Clar formula
is the {\em Clar number} of $F$, denoted by $c(F)$. In Clar's model
\cite{C}, a Clar formula is designated by depicting circles within
their hexagons (see Figure \ref{fig1-1}).

The Clar number is originally defined for benzenoid systems based on
the Clar sextet theory \cite{C} and related to Randi\'c conjugated
circuit model \cite{R1}. It is effective to measure the molecule
stability of benzenoid hydrocarbons. For two isomeric benzenoid
hydrocarbons, the one with larger Clar number is more stable. Clar
numbers of benzenoid hydrocarbons have been investigated and
computed in many papers \cite{HZ1,HZ2,KZG,ZL,Z1,Z2,Z3}. Hansen and
Zheng \cite{HZ2} introduced an integer linear program to compute the
Clar number of benzenoid hydrocarbons. Abeledo and Atkinson
\cite{AA} showed that relaxing the integer-restrictions in such a
program always yields an integral solution.

Up to now there has been no an effective method to compute  Clar
numbers of fullerene graphs. The Clar polynomial and sextet
polynomial of $\text{C}_{60}$ for counting Clar structures and
sextet patterns respectively were computed in \cite{SLZ}. This
implies that $\text{C}_{60}$ has 5 Clar formulas and Clar number 8
\cite{EB}. In addition, $\text{C}_{60}$ has a Fries structure
\cite{Fr}, i.e. a Kekul\'e structure of $\text{C}_{60}$ which avoids
double bonds in pentagons and has the possibly maximal number of
conjugated hexagons ($n/3$). Fullerene graphs with a Fries structure
are equivalent to leapfrog fullerenes or Clar type fullerenes
\cite{FP, LK}. The latter means that they have a set of disjoint
faces including all vertices, an extension of a fully-Clar
structure. Some relationships among the Clar number, the maximum
face independent number and Fries number are presented by Graver
\cite{G}. A lower bound for the Clar numbers of leapfrog fullerenes
with icosahedral symmetry was also given in \cite{G}. The same
authors of this paper \cite{ZY} showed that the Clar number of a
fullerene graph with $n$ vertices is no more than $\lfloor\frac
{n-12} 6\rfloor$, for which equality holds for infinitely many
fullerene graphs, including $\text{C}_{60}$ and $\text{C}_{70}$. We
would like to mention here that a recent paper of Kardo\v{s} et al.
\cite{KKMS} obtained a exponentially bound of perfect matching
numbers of fullerene graphs. In fact, they applied Four-Color
Theorem to show that a fullerene graph with $n\ge 380$ vertices has
a sextet pattern with at least $\frac {n-380} {61}$ hexagons.

A fullerene graph $F_n$ is {\em extremal} if its Clar number
$c(F_n)=\frac {n-12} 6$. In this paper, we characterize the extremal
fullerene graphs with at least 60 vertices (Section 3). According to
the characterization, we construct all  18 extremal fullerene graphs
with 60 vertices, including $\text{C}_{60}$ (Section 4). Our result
can show that a combination of Clar number and Kekul\'e count works
well in predicting the stability of $\text{C}_{60}$.

%This result can be  study is an extended research of \cite{ZYL}
%which shows the stability of the fullerene isomers of
%$\text{C}_{60}$ is coincident with the bi-index consisting of both
%Clar number and perfect matching count. \cite{ZYL}

\section{Definitions and Terminologies}

Let $G$ be a plane graph with vertex-set $V(G)$ and edge-set $E(G)$.
Let $|G|=|V(G)|$. For a 2-connected plane graph, every face is
bounded by a cycle. For convenience, a face is represented by its
boundary if unconfused. The boundary of the infinite face of $G$ is
also called the boundary of $G$, denoted by $\partial G$. A graph
$G$ is {\em cyclically $k$-edge-connected} if deleting less than $k$
edges from $G$ cannot separate it into two components such that each
of them contains at least one cycle. The {\em cyclic
edge-connectivity} of graph $G$, denoted by $c\lambda(G)$, is the
maximum integer $k$ such that $G$ is cyclically $k$-edge-connected.

\begin{lem}\label{lem2-1}{ \upshape{(\cite{D,QZ})}}
Let $F$ be a fullerene graph. Then $c\lambda(F)=5$. \qed
\end{lem}

From now on, let $F$ be a fullerene graph. Let $C$ be a cycle of
$F$. Lemma \ref{lem2-1} implies that the size of $C$ is larger than
4. The subgraph consisting of $C$ together with its interior is
called a {\em fragment}. A {\em pentagonal fragment} is a fragment
with only pentagonal inner face. For a fragment $B$, all 2-degree
vertices of $B$ lie on its boundary.

\begin{figure}[!hbtp]\refstepcounter{figure}\label{fig2-1}
\begin{center}
\includegraphics{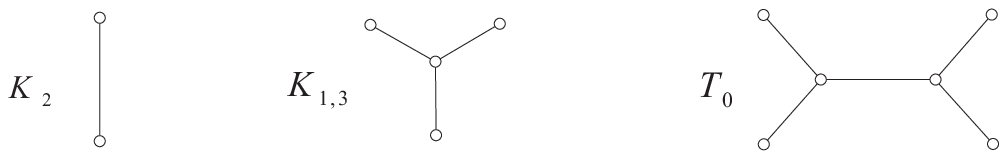}\\
{Figure \ref{fig2-1}: Trees: $K_2$, $K_{1,3}$ and $T_0$.}
\end{center}
\end{figure}

\begin{lem}\label{lem2-2}
Let $B$ be a fragment of a fullerene graph $F$ and let $W$ be the
set of all 2-degree vertices of $B$. If $0<|W|\le 4$, then
$T:=F-(V(B)\setminus W)$ is a forest and,\\
{\upshape{(1)}} $T$ is $K_2$ if $|W|=2$;\\
{\upshape{(2)}} $T$ is $K_{1,3}$ if $|W|=3$;\\
{\upshape{(3)}} $T$ is the union of two $K_2$'s, or a $3$-length
path, or $T_0$ as shown in Figure \ref{fig2-1} if $|W|=4$.
\end{lem}

\begin{proof} Since $B$ is a fragment, $\partial B$ is a cycle. For every
vertex $w\in W$, let $ww_1,ww_2\in E(\partial B)$. The neighbor of
$w$ distinct from $w_1$ and $w_2$ belongs to either $W$ or $V(F-B)$.

If $V(F)=V(B)$, then every vertex in $W$ is adjacent to exactly one
2-degree vertices in $W$. Therefore $|W|=2$ or $|W|=4$. If $|W|=2$,
then the two vertices in $W$ are adjacent. Further $T$ is a $K_2$.
If $|W|=4$, then $F$ has two more edges than $B$. If the no vertices
in $W$ are adjacent in $B$, then the two edges are disjoint and
hence $T$ is a union of two $K_2$. If there are vertices are
adjacent, then there are exactly one pair of 2-degree vertices from
$W$ adjacent since $F$ contains no 4-length cycle. It follows that
$T$ must be a 3-length path consisting of two edges in $E(F)-E(B)$
and one edge in $E(B)$.

So suppose $V(F)\setminus V(B)\ne \emptyset$. Let $S$ be a set of
the edges joining the vertices in $W$ and their neighbors in $F-B$.
Since every vertex in $W$ has at most one neighbor in $F-B$, we have
$|S|\le |W|$. So $S$ separates $B$ from $F-B$. By Lemma
\ref{lem2-1}, $F-B$ has no cycles since $|W|\le 4$.

Suppose to the contrary that $T$ has at least one cycle $C$. Then
$C\cap \partial B\ne \emptyset$ since $F-B$ is a forest. We draw $F$
on the plane such that $B$ lies outside of $C$. Then $C$ together
with its interior is a subgraph of $T$. We may thus assume that $C$
bounds a face of $F$ within $T$. Since $F$ is cubic, every component
of $C\cap
\partial B$ is an edge joining two vertices in $W$. By $0<|W|\le 4$,
$C\cap
\partial B$ has at most two components.

If $C\cap \partial B$ has two components, then $|W|=4$ and $C$
contains all vertices in $W$. Let $w_1,w_2,w_3,w_4$ be the four
vertices in $W$ and let $w_1w_2$ and $w_3w_4$ be the two components
of $C\cap B$. Let $w'_i$ be another neighbor of $w_i$ on $\partial
B$. %Further, $w'_i\ne w_j$ for $i,j\in\{1,2,3,4\}$.
Then $\{w_iw'_i|i=1,2,3,4\}$ separates $C$ from $F-C$, contradicting
Lemma \ref{lem2-1}.

So suppose $C\cap \partial B$ has only one component $w_1w_2$. Then
$C-\{w_1,w_2\}$ is a path in a component $T_1$ of $F-B$. Further
$T_1$ has at least $|C|-2$ vertices. If $T_1$ has a 3-degree vertex,
then it has at least three leaves. Since every leaf of $T_1$ is
adjacent to two vertices in $W$, we have $|W|\ge 6$ which
contradicts that $|W|\le 4$. So $T_1$ is a path. Then $T_1$ has at
least $|C|-4$ 2-degree vertices. Hence vertices in $V(T_1)$ have at
least $4+|C|-4$ neighbors in $W$. So $|W|\ge 4+|C|-4\ge 5$ which
also contradicts that $|W|\le 4$. So $T$ is a forest.

Let $l$ and $x$ be the number of leaves and the number of components
of $T$, respectively. Then $l=|W|\le 4$. Since $F$ is cubic,
$2(|T|-x)=3(|T|-l)+l$. Then $l-2x=|T|-l> 0$ since $F-B\ne \emptyset$
and $W\ne \emptyset$. Hence $4\ge l> 2x\ge 2$. So we have $x=1$.
Hence $T$ is a tree. So if $l=3$, then $T$ is $K_{1,3}$. If $l=4$,
then $|T|-l=2$. Hence $T$ is isomorphic to $T_0$.\end{proof}

For a face $f$ of a connected plane graph, its boundary is a closed
walk. For convenience, a face $f$ is often represented by its
boundary if unconfused. Note that a pentagon or a hexagon of a
fullerene graph $F$ must bound a face since $F$ is cyclic
5-edge-connected \cite{D,ZY}. Let $G$ be a subgraph of a fullerene
graph $F$. A face $f$ of $F$ {\em adjoins} $G$ if $f$ is not a face
of $G$ and $f$ has at least one edge in common with $G$. Now suppose
$G$ has no 1-degree vertices. Let $f'$ be a face of $G$ with
2-degree vertices on its boundary. Since $F$ is cubic and
3-connected, $f'$ has at least two 2-degree vertices. A path $P$ on
the boundary of $f'$ connecting two 2-degree vertices is {\em
degree-saturated} if $P$ contains no 2-degree vertices of $G$ as
intermediate vertices. Since every face of $F$ has a size of at most
six, the length of $P$ is no more than five.

\begin{pro}\label{pro2-3}
Let $G$ be a subgraph of a fullerene graph $F$. Let $f$ be a face of
$G$ with 2-degree vertices and $P$ be a degree-saturated path of $G$
on the boundary of $f$. Then the length of $P$ is no more than 5.
\end{pro}

\begin{figure}[!hbtp]\refstepcounter{figure}\label{fig2-2}
\begin{center}
\scalebox{0.9}{\includegraphics{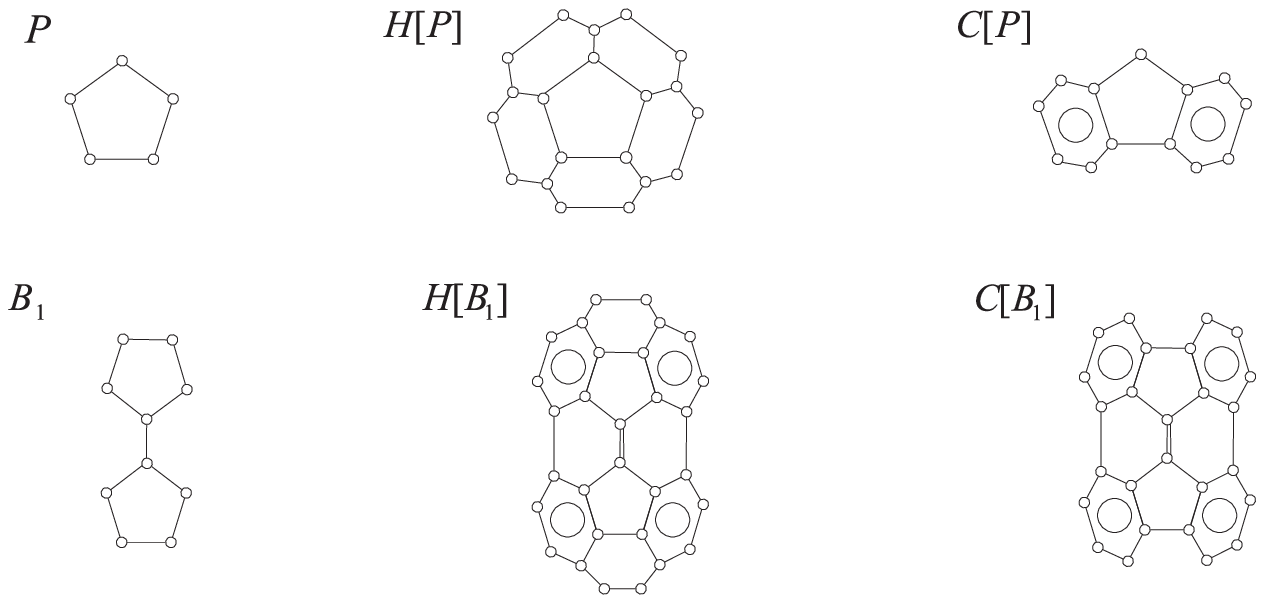}}\\
{Figure \ref{fig2-2}: The hexagon extensions and Clar extensions of
$P$ and $B_1$.}
\end{center}
\end{figure}

Let $f_1,f_2,...,f_k$ be the faces of $F$ adjoining $G$. The
subgraph $T[G]:=G\cup (\cup_{i=1}^k f_i)$ is called the {\em
territory} of $G$ in $F$. If for every $i\in \{1,2,...,k\}$, the
face $f_i$ ($i=1,...,k$) is a hexagon, the territory is also called
a {\em hexagon extension} of $G$ and is denoted by $H[G]$ (see
Figure \ref{fig2-2}). A subgraph $G$ is {\em maximal} in $F$ if
$H[G]\subset F$. We are particularly interested in the maximal
pentagonal fragments. Denote the number of 2-degree vertices of $G$
by $w(G)$. Let $B$ and $B'$ be two fragments such that $w(B)\ge
w(B')$. Let $P$ and $P'$ be two degree-saturated paths of $\partial
B$ and $\partial B'$, respectively. Suppose $|P|\le |P'|$. Let $f$
and $f'$ be two faces adjoining $B$ and $B'$ along $P$ and $P'$,
respectively. It is readily seen that $w(B\cup f)\ge w(B'\cup f')$
if $|f|\ge |f'|$. Applying this argument for the territory $T[B]$
and the hexagon extension $H[B]$ of $B$, we immediately have the
following proposition.

\begin{pro}\label{pro2-4}
Let $B$ be a fragment of a fullerene graph $F$ and let $T[B]$ and
$H[B]$ be the territory and the hexagon extension of $B$,
respectively. Then $w(T[B])\le w(H[B])$.
\end{pro}

A subgraph (or a set of vertices) $S$ of $F$ {\em meets} a subgraph
$G$ of $F$ if $S\cap G\ne \emptyset$. Let $G-S$ be the subgraph
obtained from $G$ by deleting all vertices in $S$ together with all
edges incident with them. Let $H[G]$ be the hexagon extension of $G$
and $\mathcal H$ be a set of mutually disjoint hexagons of $H[G]$.
Let
$$\mathscr
S(G):=\{\mathcal H|\ G-\mathcal H \mbox{ has a matching which covers
all remaining 3-degree vertices of } G\}.$$ For any $\mathcal H\in
\mathscr S(G)$, let $U_{\mathcal H}(G):=V(G)\setminus V(\mathcal
H)$. An $\mathcal H\in \mathscr S(G)$ is called a {\em Clar set} of
$H[G]$ if $|U_{\mathcal H}(G)|\le |U_{\mathcal H'}(G)|$ for all
$\mathcal H'\in \mathscr S(G)$. A Clar set $\mathcal H$ of $H[G]$ is
{\em normal} if $G-\mathcal H$ has a perfect matching. For a
fullerene graph $F$, its hexagon extension is itself and a Clar
formula of $F$ is a normal Clar set. (See Figure \ref{fig2-2}: the
hexagons in Clar sets of the hexagon extensions of a pentagon $P$
and $B_1$ are depicted by circles; the Clar set of $H[B_1]$ is
normal.)

\begin{defi}{\rm Let $G\subseteq F$ and $\mathcal H$ be a Clar set
of $H[G]$. A {\em Clar extension} $C[G]$ of $G$ is the subgraph
induced by $V(\mathcal H)\cup V(G)$. A Clar extension $C[G]$ is {\em
normal} if $\mathcal H$ is normal.}
\end{defi}

The Clar extensions of $P$ and $B_1$ are illustrated in Figure
\ref{fig2-2}. The Clar extension of a fullerene graph $F$ is itself.
Let $\mathcal H$ be a Clar formula of $F$ and $U_{\mathcal
H}:=V(F)\setminus V(\mathcal H)$. The following result is from
\cite{ZY}.

\begin{lem} {\rm (\cite{ZY}, Lemma 2)}\label{lem2-6}
If a subgraph $G$ of a fullerene graph $F$ has at least $k$
pentagons, then $|V(G)\cap U_{\mathcal H}|\ge k$.\qed
\end{lem}

Lemma \ref{lem2-6} can be generalized as the following result.

\begin{lem}\label{lem2-7}
Let $G$ be a subgraph of a fullerene graph $F$ with $k$ pentagons
and $\mathcal H$ be a Clar set of $H[G]$. Then $|U_{\mathcal H}
(G)|\ge k$.
\end{lem}

\begin{proof} Let $G$ be a subgraph of $F$ with $k$ pentagons and $\mathcal
H$ be a Clar set of $H[G]$. We proceed by induction on $k$. If
$k=1$, $|U_{\mathcal H}(G)|\ge 1$ since every pentagon has at least
one vertex not in $\mathcal H$. So suppose the conclusion holds for
smaller $k$.

If $G$ has a 2-degree vertex $v$ in $U_{\mathcal H}(G)$, then $G-v$
has at least $k-1$ pentagons. By inductive hypothesis, $|U_{\mathcal
H}(G-v)|\ge k-1$. So $|U_{\mathcal H}(G)|=|U_{\mathcal H}(G-v)|+1\ge
k$ and the lemma holds.

So suppose all vertices in $U_{\mathcal H}(G)$ are 3-degree vertices
of $G$; that is, $G-V(\mathcal H)$ has a perfect matching. The proof
of this case follows directly from the proof of Lemma \ref{lem2-6}
(Lemma 2 in \cite{ZY}). \end{proof}

A subgraph $G$ with $k$ pentagons is {\em extremal} if $|U_{\mathcal
H}(G)|=k$ where $\mathcal H$ is a Clar set of $H[G]$. Both $P$ and
$B_1$ are extremal (see Figure \ref{fig2-2}). Note that the every
subgraph induced by pentagons of an extremal fullerene graph must be
extremal. Hence extremal subgraphs play a key role in characterizing
extremal fullerene graphs.

\section{Extremal fullerene graphs}

In this section, we are going to characterize extremal subgraphs
induced by pentagons of fullerene graphs and finally establish a
characterization of the extremal fullerene graphs with at least 60
vertices.

\begin{figure}[!hbtp]\refstepcounter{figure}\label{fig3-1}
\begin{center}
\includegraphics{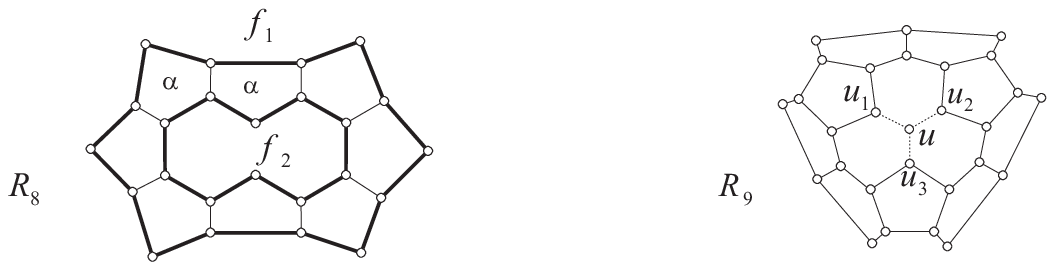}\\
{Figure \ref{fig3-1}: Pentagonal rings: $R_8$ and $R_9$.}
\end{center}
\end{figure}

From now on, let $F_n$ be a fullerene graph with $n$ vertices. %Let
%$G$ be a subgraph of $F_n$. A path $P$ on a face $f$ of $G$
%connecting two 2-degree vertices is {\em saturated} if every inner
%vertex of $P$ is a 3-degree vertex of $G$.
A {\em pentagonal ring} $R_k$ is a subgraph of $F_n$ consisting of
$k$ pentagons $P_0,P_1,...,P_{k-1}$ such that $P_i\cap P_j\ne
\emptyset$ if and only if $|i-j|=1$ where $i,j\in \mathbb Z_k$ (see
Figure
\ref{fig3-1}). % any two consecutive pentagons (mod $k$) have a common
%edge and any three pentagons have no common vertices.
Since $F_n$ has exactly 12 pentagons and $c\lambda(F_n)=5$, we
deduce that $5\le k\le 12$.

\begin{lem}\label{lem3-1}
If $F_n$ contains a pentagonal ring $R_k$ with $7\le k\le 12$, then
$n\leq 52$.
\end{lem}
\begin{proof}  Let $R_k\subset F_n$ be a pentagonal ring. Let $f_1,
f_2\notin \{P_0,...,P_{k-1}\}$ be two faces of $R_k$. We may assume
that $f_1$ is the infinite face. For $i\in\{1,2\}$, let $x_i$ be the
numbers of 2-degree vertices on the boundary of $ f_i$.

Let $B$ be the fragment consisting of $f_1$ together with its
interior. Let $m_5$ and $m_6$ be the number of pentagons and
hexagons of $B$, respectively. By Euler's formula,
\[\nu-e+m_5+m_6=1
\]
where $\nu, e$ are the vertex number and the edge number of $B$,
respectively. On the other hand,
\[
2x_1+3(\nu-x_1)=2e=5m_5+6m_6+2x_1+x_2.
\]
Hence, $m_5=6+x_2$. Since $m_5\ge k=x_1+x_2$, it follows that
$x_1\le 6$. Since $7\le k\le 12$ and $F$ is 3-connected, $x_2\ge 2$.
It can be verified that $H[B]$ has at most four 2-degree vertices on
its boundary. Let $B'$ be the fragment consisting of $f_2$ together
with its interior. By Lemma \ref{lem2-2} and Proposition
\ref{pro2-4}, $|V(F-B')|\le 6\times 6-2\times 6+2=26$ since there
are at most six faces adjoining $B$ and any two adjacent faces share
at least one edge. A similar discussion results in $|V(B')|\le 26$
since $F$ can be drawn on the plane such that $f_2$ is the infinite
face of $R_k$. So $\nu\le |V(F-B')|+|V(B')|\le 52$. \end{proof}

\begin{figure}[!hbtp]\refstepcounter{figure}\label{fig3-2}
\begin{center}
\includegraphics{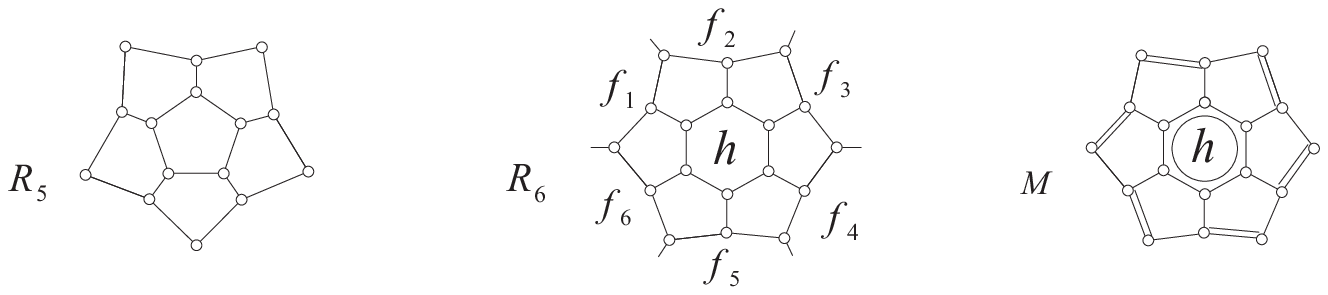}\\
{Figure \ref{fig3-2}: Pentagonal rings $R_5$, $R_6$ and a matching
$M$ of $R_6$.}
\end{center}
\end{figure}

The following observations show that a subgraph $G$ of $F_n$ (except
$F_{24}$) is not extremal if it contains $R_5$ and $R_6$ as
subgraphs. Recall that the territory and the hexagon extension of
$G$ is denoted by $T[G]$ and $H[G]$, respectively. For a Clar set
$\mathcal H$ of $H[G]$, define $U_{\mathcal H}(G):=V(G)\setminus
V(\mathcal H)$. Let $R_5$ and $R_6$ be the pentagonal rings depicted
in Figure \ref{fig3-2}.

\begin{obs}\label{obs1}
Let $\mathcal H$ be a Clar set of $H[R_5]$. Then $|U_{\mathcal
H}(R_5)|\ge 12$.
\end{obs}

The proof of Observation \ref{obs1} is omitted here since it is
similar to the proof of Lemma 1 of \cite{ZY}.

\begin{obs}\label{lem3-3}
Let $G$ be a subgraph of a fullerene graph. If $R_5\subseteq G$,
then $G$ is not extremal.
\end{obs}

\begin{proof} Suppose to the contrary that $G$ is extremal . Since
$R_5\subseteq G$, we have that $|U_{\mathcal H}(G)|\ge |U_{\mathcal
H}(R_5)|\ge 12$ by Observation \ref{obs1}. Hence $G$ has 12
pentagons since $G$ is extremal. Clearly, every pentagon contains at
least one vertex in $U_{\mathcal H}(G)$ and at least one pentagon
does not adjoin $R_5$. So $|U_{\mathcal H}(G)|\ge |U_{\mathcal
H}(R_5)|+1=13$ which contradicts that $G$ is extremal . \end{proof}

\begin{obs}\label{lem3-4}
Let $G$ be a subgraph of a fullerene graph $F_n$ with $n\ne 24$. If
$R_6\subseteq G$, then $G$ is not extremal .
\end{obs}

\begin{proof} Let $\mathcal H$ be a Clar set of the hexagon extension $H[G]$
of $G$. Enumerate clockwise the six faces of $F_n$ adjoining $R_6$
as $f_1,...,f_6$ (see Figure \ref{fig3-2}). Since $G\subseteq F_n$
($n\ne 24$), not all $f_i$ ($1\le i\le 6$) are pentagons. Let
$r:=|\mathcal H \cap \{f_1,\cdots, f_6\}|$ and $h$ be the central
hexagon of $R_6$.

If $h\notin \mathcal H$, then $|U_{\mathcal H}(G)|\ge |U_{\mathcal
H}(R_6)|=18-3r$ since every $f_i$ contains three vertices in
$V(R_6)$. If $r\le 1$, then $G$ is not extremal  since $|U_{\mathcal
H}(G)|\ge 15$ and $G$ has at most 12 pentagons. So $2\le r \le 3$.
If $r=3$, say $f_1, f_3, f_5\in \mathcal H$, then $R_6- \mathcal H$
has no matchings which cover all remaining 3-degree vertices of
$R_6$, contradicting that $\mathcal H$ is a Clar set. So suppose $r
= 2$. Then $G$ has exact 12 pentagons. Over these 12 pentagons, at
least two pentagons do not adjoin $R_6$. Since every pentagon
contains at least one vertex in $U_{\mathcal H}(G)$, it holds that
$|U_{\mathcal H}(G)|\ge 12+1=13$. Hence $G$ is not extremal .

So suppose $h\in \mathcal H$. Then all 3-degree vertices on
$\partial R_6$ of $R_6$ have to match all 2-degree vertices on
$\partial R_6$ in $G-\mathcal H$ (see Figure \ref{fig3-2}, $R_6$
with a matching $M$). So $|U_{\mathcal H}(G)|\ge |V(\partial
R_6)|=12$. So suppose $G$ has 12 pentagons. Since $G\subseteq F_n\ne
F_{24}$, at least one pentagon in $G$ does not adjoin $R_6$ and has
at least one vertex in $U_{\mathcal H}(G)$. Immediately,
$|U_{\mathcal H}(G)|\ge 12+1=13$. So $G$ is not extremal .
\end{proof}

By the above observations and Lemma \ref{lem3-1}, an extremal
fullerene graph with at least 60 vertices does not contain a
pentagonal ring as a subgraph. If a connected component of the
subgraph induced by pentagons of $F_n$ with $n\ge 60$ is extremal,
then it must be a pentagonal fragment.

\begin{figure}[!hbtp]\refstepcounter{figure}\label{fig3-3}
\begin{center}
\scalebox{0.92}{\includegraphics{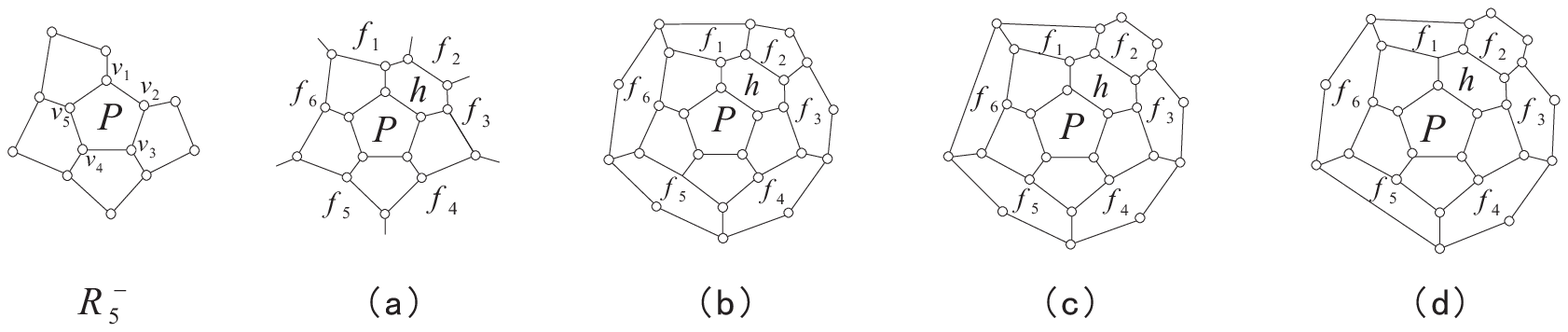}}\\
{Figure \ref{fig3-3}: The fragment $R_5^-$ and illustration for the
proof of Lemma \ref{lem3-5}.}
\end{center}
\end{figure}

Let $R_5^-$ be the pentagonal fragment arising from $R_5$ by
deleting one 2-degree vertex together with two edges incident with
it (see Figure \ref{fig3-3}).

\begin{lem}\label{lem3-5}
Let $G$ be a subgraph of $F_n$ with $n\ge 40$. If $R_5^-\subseteq
G$, then $G$ is not extremal .
\end{lem}

\begin{proof} Let $G\subseteq F_n$ with $k$ pentagons and $\mathcal H$ be a
Clar set of $H[G]$. Suppose to the contrary that $G$ is extremal .
By Lemma \ref{lem3-3} and \ref{lem3-4}, $R_5\nsubseteq G$ and
$R_6\nsubseteq G$. Let $P:=v_1v_2...v_5v_1$ be the pentagon of
$R_5^-$ meeting all other pentagons of $R_5^-$ as shown in Figure
\ref{fig3-2}. Let $h$ be the hexagon of $F_n$ adjoining $R_5^-$
along $v_1v_2$ since $R_5\nsubseteq G$. Let $f_1,f_2,...,f_6$ be the
faces of $F_n$ adjoining $R_5^-\cup h$ as shown in Figure
\ref{fig3-3} (a).

Let $r$ be the number of pentagons in $\{f_1,...,f_6\}$ and
$H[R_5^-\cup h]$ be the hexagon extension of $R_5^-\cup h$. Clearly,
$H[R_5^-\cup h]$ has seven 2-degree vertices. If $r\ge 3$, then the
territory $T[R_5^-\cup h]$ of $R_5^-\cup h$ has at most four
2-degree vertices on its boundary. By Lemma \ref{lem2-2}, $n\le
|V(T[R_5^-\cup h])|+2\le 26+2=28$, contradicting that $n\ge 40$. So
suppose $r\le 2$.

If $r=2$, then the boundary of $T[R_5^-\cup h]$ has five 2-degree
vertices which separate $\partial (T[R_5^-\cup h])$ into five
degree-saturated paths. If $f_2$ is a pentagon, $\partial
(T[R_5^-\cup h])$ has four 2-length degree-saturated paths and one
3-length degree-saturated path (see Figure \ref{fig3-3} (b)). Then
the hexagon extension $H[T[R_5^-\cup h]]$ has only four 2-degree
vertices on its boundary. By Lemma \ref{lem2-2} and Proposition
\ref{pro2-4}, $n\le |V(H[T[R_5^-\cup h]])|+2=|V(T[R_5^-\cup
h])|+9+2=38$, contradicting $n\ge 40$. So suppose $f_2$ is a
hexagon. If the two pentagons in $\{f_1,f_3,f_4,f_5,f_6\}$ are
adjacent, $\partial (T[R_5^-\cup h])$ has one 1-length
degree-saturated path and three 2-length degree-saturated paths and
one 4-length degree-saturated path (see Figure \ref{fig3-3} (c)).
Then $H[T[R_5^-\cup h]]$ has 35 vertices and three 2-degree
vertices. By Lemma \ref{lem2-2} and Proposition \ref{pro2-4}, $n\le
|V(H[T[R_5^-\cup h]])|+1=36$, contradicting that $n\ge 40$. So we
suppose the two pentagons in $\{f_1,f_3,f_4,f_5,f_6\}$ are not
adjacent. Then $\partial (T[R_5^-\cup h])$ has one 1-length
degree-saturated path and two 2-length degree-saturated paths and
two 3-length degree-saturated paths (see Figure \ref{fig3-3} (d)).
Further, $H[T[R_5^-\cup h]]$ has 36 vertices and four 2-degree
vertices. By Lemma \ref{lem2-2} and Proposition \ref{pro2-4}, $n\le
|V(H[T[R_5^-\cup h]])|+2=38$, also contradicting that $n\ge 40$.

So $r=1$. Suppose $h\in \mathcal H$. Then $f_1,f_2,f_3\notin
\mathcal H$. If $f_4,f_6\in \mathcal H$, then $R_5^-\cup h-\mathcal
H$ has no matchings which cover all remaining 3-degree vertices of
$R_5^-\cup h$, a contradiction. So at most one of $\{f_4,f_5,f_6\}$
belongs to $\mathcal H$. Hence, $|U_{\mathcal
H}(R_5^-)|=|V(R_5^-\cup h)|-|V(\mathcal H)|\ge 16-9=7$. Since $G$
has $k$ pentagons and $r=1$, it holds that $G-R_5^-$ has at least
$k-6$ pentagons. By Lemma \ref{lem2-7}, $|U_{\mathcal
H}(G-R_5^-)|\ge k-6$. Hence $| U_{\mathcal H}(G)|=|U_{\mathcal
H}(R_5^-)|+|U_{\mathcal H}(G-R_5^-)|\ge 7+k-6=k+1$. Hence $G$ is not
extremal.

Now suppose that $h\notin \mathcal H$. Since both $(R_5^-\cup
h)-(\cup_{i=1,3,5} f_i)$ and $(R_5^-\cup h)-(\cup_{i=2,4,6} f_i)$
have no perfect matchings, at most two faces of $f_1,...,f_6$ belong
to $\mathcal H$. So $|U_{\mathcal H}(R_5^-)|\ge 14-6=8$. On the
other hand, $G-R_5^-$ has $k-6$ pentagons since $r=1$. Hence, by
Lemma \ref{lem2-7}, $|U_{\mathcal H}(G)|=|U_{\mathcal
H}(R_5^-)|+|U_{\mathcal H}(G-R_5^-)|\ge 8+ k-6\ge k+2$, a
contradiction. Hence $G$ is not extremal.\end{proof}

For a pentagonal fragment $B$, let $\gamma (B)$ be the minimum
number of pentagons adjoining a common pentagon in $B$. Let $B^*$ be
the inner dual of $B$. Then $\gamma (B)$ is the the minimum degree
of $B^*$. For example, $\gamma (R_5)=3$ and $\gamma(R_5^-)=2$.

\begin{lem}\label{lem3-6}
Let $B$ be a pentagonal fragment of a fullerene graph $F$. Then:\\
{\upshape(1)} $R_5\subseteq B$ if $\gamma(B)\ge 3$;\\
{\upshape(2)} $B$ has a pentagon adjoining exactly two adjacent
pentagons of $B$ if $\gamma(B)=2$.
\end{lem}

\begin{proof} Let $B^*$ be the inner dual of $B$. Then $B^*$ is a simple
connected graph and every inner face of $B^*$ is a triangle. Let
$\delta(B^*)$ be the minimum degree of $B^*$. Then
$\delta(B^*)=\gamma(B)$.

Suppose to the contrary that $R_5\nsubseteq F$; that is, $B^*$ is an
outer plane graph. It suffices to prove that $\delta(B^*)\le 2$ and
$B^*$ has a 2-degree vertex on a triangle of $B^*$ if
$\delta(B^*)=2$. If $\delta(B^*)=1$, the assertion already holds. So
suppose $\delta(B^*)=2$. Let $G$ be a maximal 2-connected subgraph
of $B^*$ such that $G$ is connected to $F-G$ by an edge $e$. If
$B^*$ is 2-connected, let $G=B^*$. Then every inner face of $G$ is a
triangle. So it suffices to prove that $G$ has two 2-degree
vertices.

Let $C$ be the boundary of $G$. Let $v_0,v_1,v_2,...,v_{n-1}$ be all
vertices of $G$ appearing clockwise on $C$. If $n=3$, then $G$ is a
triangle and the assertion is true. So suppose $n>3$. Since every
inner face of $G$ is a triangle, then $G$ has 3-degree vertices.
Without loss of generality, let $v_0$ be a 3-degree vertex such that
$v_0v_k$ is a chordal of $C$ where $k\ne 1,n-1$. Let $v_jv_{j'}$ be
a chordal of $C$ such that $k\le j<j+1<j'\le n\equiv 0$ (mod $n$)
and $|j'-j|$ is minimal. Then the cycle $v_jv_{j+1}\cdots
v_{j'-1}v_{j'}v_j$ bounds an inner face. So it is a triangle and
$v_{j+1}$ is a 2-degree vertex on the triangle
$v_jv_{j+1}v_{j'}v_j$. On the other hand, let $v_{i}v_{i'}$ be a
chordal of $C$ such that $0\le i<i+1<i'\le k$ and $i'-i$ is minimal.
A similar analysis implies that $v_{i+1}$ is a 2-degree vertex on
the triangle $v_iv_{i+1}v_{i'}v_{i}$. At most one of $v_{j+1}$ and
$v_{i+1}$ is an end of the edge $e$ joining $G$ to $F-G$. So $B$ has
a 2-degree vertex on a triangle of $B$. This completes the proof of
the lemma. \end{proof}

\begin{lem}\label{lem3-7}
Let $B$ be a pentagonal fragment with $\gamma (B)\ge 2$. Then $B$ is
not extremal .
\end{lem}
\begin{proof} Let $k$ be the number of pentagons of $B$ and $\mathcal H$ be
a Clar set of $H[G]$. Use induction on $k$ to prove it. The minimum
pentagonal fragment $B_0$ with $\gamma(B_0)\ge 2$ consists of three
pentagons such that they adjoin each other. It is easy to verify
that $B_0$ is not extremal . So we may suppose $k\ge 4$ and the
lemma holds for smaller $k$. If $R_5\subseteq B$, then $B$ is not
extremal according to Lemma \ref{lem3-3}. By Lemma \ref{lem3-6}, we
may assume $\gamma (B)=2$ and let $p:=v_1v_2v_3v_4v_5v_1$ be a
pentagon of $B$ adjoining two pentagons $p_1$ and $p_2$ such that
$p_1\cap p=v_3v_4$ and $p_2\cap p=v_4v_5$.

Let $h_1,h_2,h_3$ be the three hexagons of $F_n$ adjoining $p$ as
illustrated in Figure \ref{fig3-4} (a). If one of $v_1$ and $v_2$
belongs to $U_{\mathcal H}(B)$, then $B':=B-\{v_1,v_2\}$ has at
least $k-1$ pentagons and $\gamma(B')\ge 2$. So $B'\notin \mathscr
B_{\ge 60}$. By inductive hypothesis, $B'$ is not extremal  and
hence $| U_{\mathcal H}(B')|\ge k$. Hence $|U_{\mathcal H}(B)|\ge
|V(B')\cap U_B|+1\ge k+1$. That means $B$ is also not extremal . So
suppose $v_1,v_2\in V(\mathcal H)$. Then either $h_2\in \mathcal H$
or $h_1, h_3\in \mathcal H$.

\begin{figure}[!hbtp]\refstepcounter{figure}\label{fig3-4}
\begin{center}
\includegraphics{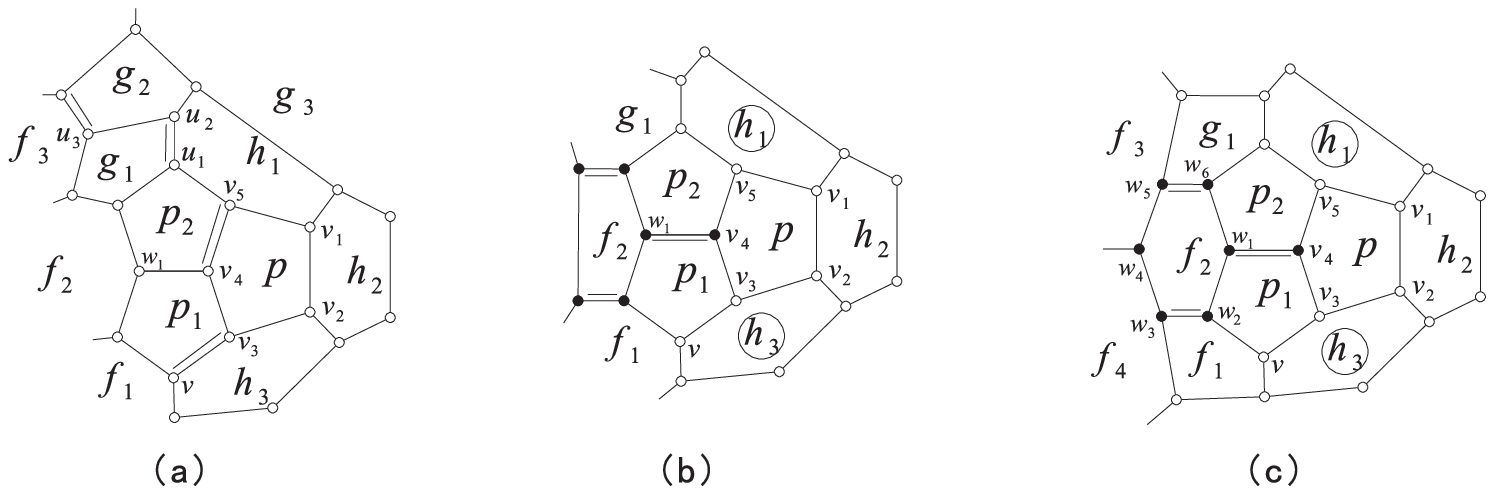}\\
{Figure \ref{fig3-4}: Illustration for the proof of Lemma
\ref{lem3-7}.}
\end{center}
\end{figure}

\noindent {\em Case 1:} $h_2\in \mathcal H$. Then $v_3,v_4,v_5\in
U_{\mathcal H}(B)$ and all of them are covered by $M_B$. Let $f_1,
f_2$ be the other two faces adjoining $p_1$ as shown in Figure
\ref{fig3-4} (a). Let $w_1v_4=p_1\cap p_2$. If $f_2\notin \mathcal
H$, then $S=\{w_1,v_3,v_4,v_5\}\subseteq U_{\mathcal H}(B)$, a
contradiction. So suppose $f_2\in \mathcal H$. So either $v_3v_4\in
M_B$ or $v_4v_5\in M_B$. By symmetry, we may assume $v_4v_5\in M_B$.
Let $vv_3=p_1\cap h_3$. Then $vv_3\in M_B$. Since
$\{v_3,v_4,v_5,v\}\subseteq U_{\mathcal H}(B)$ should meet at least
four pentagons, $f_1$ is a pentagon.

Let $g_1,g_2,g_3$ be the faces adjoining $h_1$ as illustrated in
Figure \ref{fig3-4} (a), and let $u_1u_2=g_1\cap h_1$ and
$u_2u_3=g_1\cap g_2$. Since $g_1\notin \mathcal H$, we have $u_1\in
U_{\mathcal H}(B)$. Since $\{v_3,v_4,v_5,u_1\}\subseteq U_{\mathcal
H}(B)$ meets at least four pentagons, $g_1$ is a pentagon. Hence
$u_1u_2\in M_B$. So $\{v_3,v_4,v_5,u_1,u_2\}\subseteq U_{\mathcal
H}(B)$. Further $g_2$ is also a pentagon. Let $f_3$ be the face
adjoining $g_1,g_2$ and $f_2$. Then $f_3\notin \mathcal H$ since it
is adjacent with $f_2$. Further $f_3$ is a pentagon since
$\{v_3,v_4,v_5,u_1,u_2,u_3\}\subseteq U_{\mathcal H}(B)$ meets at
least six pentagons.

Let $B':=B-(V(P)\cup \{w_1\})$. If $B'$ is connected, then the
pentagons in $B'$ connecting $f_1$ and $g_1$ together with $p_1,p_2$
form a pentagonal ring in $B$, contradicting that $B$ is a
pentagonal fragment. Let $B_1,...,B_r$ be all components of $B'$
such that $g_1\subseteq B_1$. Use $k_i$ to denote the number of
pentagons in $B_i$, then $k=\sum_{i=1}^r k_i+3$. For $B_1$, we have
$\gamma(B_1)\ge 2$ and hence $B_1\notin \mathscr B_{\ge 60}$. By
inductive hypothesis, $B_1$ is not extremal . So $|U_{\mathcal
H}(B_1)|\ge k_1+1$. By Lemma \ref{lem2-7}, $|U_{\mathcal
H}(B)|=\sum_{i=1}^r |U_{\mathcal H}(B_i)|+3\ge (k_1+1)+\sum_{i=2}^r
k_i+3=k+1$. So $B$ is not extremal .

\noindent {\em Case 2:} $h_1,h_3\in \mathcal H$. Let $w_1v_3=p_1\cap
p_2$. Then $w_1v_3\in M_B$. Let $f_1, f_2, g_1$ be the other three
faces adjoining $p_1$ or $p_2$ (see Figure \ref{fig3-4} (b)). If
$f_2$ is a pentagon, then $V(f_2)\subseteq U_{\mathcal H}(B)$ since
$f_1,g_1\notin \mathcal H$. Hence $V(f_2)$ meets at least five
pentagons. That means $f_2$ is adjacent with at least four pentagons
in $B$, forming a $R_5^-$ in $B$. So $B$ is not extremal  by Lemma
\ref{lem3-5}. So suppose $f_2$ is a hexagon. Clearly, $f_2\notin
\mathcal H$ since $w_1\in U_{\mathcal H}(B)$.

Since $\gamma (B)=2$, both $g_1$ and $f_1$ are pentagons. Let
$f_2:=w_1w_2w_3w_4w_5w_6w_1$ and let $f_3, f_4$ be the other two
faces adjoining $f_2$ (see Figure \ref{fig3-4} (c)). Since $B$ is a
pentagonal fragment, at most one of $f_3$ and $f_4$ is a pentagon.
If exactly one of them is a pentagon, then $V(f_2)\subseteq
U_{\mathcal H}(B)$ meets only five pentagons, a contradiction. So
suppose both of them are hexagons. Then
$\{w_1,w_2,w_3,w_5,w_6\}\subseteq U_{\mathcal H}(B)$ meets only four
pentagons, also a contradiction. So $B$ is not extremal.\end{proof}

\begin{figure}[!hbtp]\refstepcounter{figure}\label{fig3-5}
\begin{center}
\includegraphics{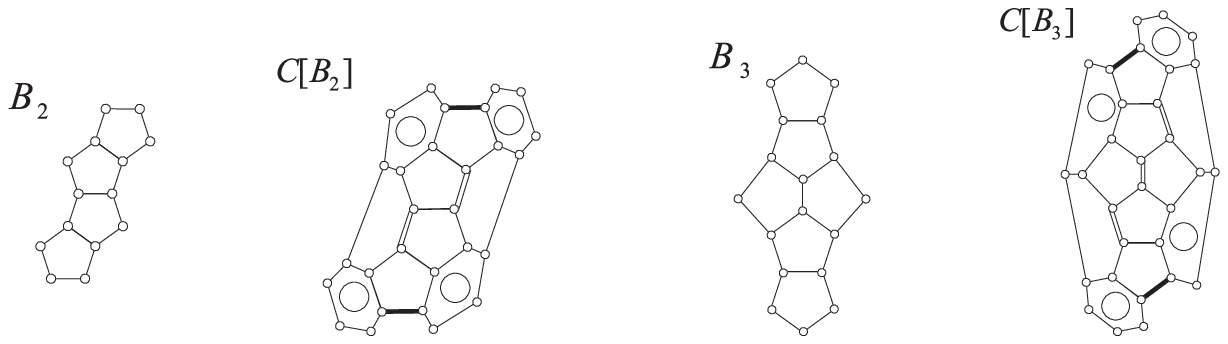}\\
{Figure \ref{fig3-5}: Extremal pentagonal fragments $B_{2}, B_{3}$
and their Clar extensions $C[B_2], C[B_3]$.}
\end{center}
\end{figure}

Now, we are going to characterize extremal pentagonal fragments. Let
$B_2, B_3$ be the two pentagonal fragments illustrated in Figure
\ref{fig3-5}. Clearly, the Clar extensions of $B_{2}$ and $B_{3}$
are normal. It is easy to see that $P, B_{2}$ and $B_{3}$ are
extremal. Up to isomorphism, $C[P], C[B_2]$ and $C[B_3]$ are unique.

\begin{figure}[!hbtp]\refstepcounter{figure}\label{fig3-6}
\begin{center}
\includegraphics{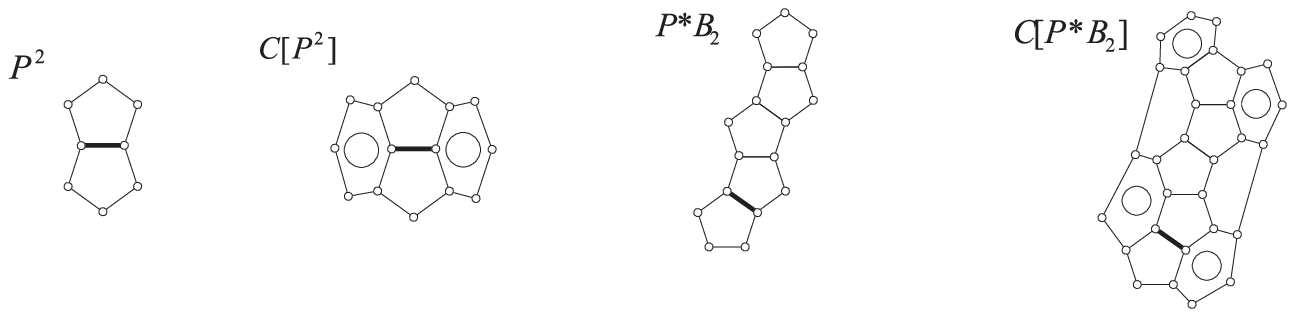}\\
{Figure \ref{fig3-6}: The pasting operation: $P^2$, $P\ast B_{2}$
and their Clar extensions.}
\end{center}
\end{figure}

Let $G_1,G_2,G_3$ and $G_4$ be graphs. We say that $G_1$ arises from
{\em pasting} $G_2$ and $G_3$ along $G_4$ if  $G_1=G_2\cup G_3$ and
$G_4=G_2\cap G_3$. Let $B$ be a fragment isomorphic to one of $P,
B_{2}$ and $B_{3}$ and let $C[B]$ be a Clar extension of $B$. An
edge of $B$ is called {\em pasting edge} if it lies on the boundary
of $C[B]$ and two end-vertices belong to $V(\mathcal H)$ where
$\mathcal H$ is the Clar set of $C[B]$. The thick edges of $P,
B_{2}$ and $B_{3}$ illustrated in Figure \ref{fig3-5} and Figure
\ref{fig3-6} are pasting edges. We can paste $P, B_{2}, B_{3}$ with
each other or itself along the pasting edges to form a new
pentagonal fragment. Use ``$\ast$" to denote the pasting operation.
Up to isomorphism, $P\ast P$ and $P\ast B_{2}$ are illustrated in
Figure \ref{fig3-6}. Simply, use $X^k$ to denote the graph obtained
pasting $k$ graphs isomorphic to $X$ along the pasting edges
together, where $X\in\{P, B_2,B_3\}$. Note that the pasting
operation does not always yield a subgraph of a fullerene graph. Let
$\mathscr B$ be the set of all maximal pentagonal fragments, which
are subgraphs of some fullerene graph, generated from the pasting
operation. Let $\mathscr B_{\ge 60}\subset \mathscr B$ such that
$B\subset F_n$ ($n\ge 60$) for any $B\in \mathscr B_{\ge 60}$.

\begin{lem}\label{lem3-8}
$B_{2}\ast B_{3}, B_{3}^2\notin \mathscr B$.
\end{lem}

\begin{figure}[!hbt]\refstepcounter{figure}\label{fig3-7}
\begin{center}
\includegraphics{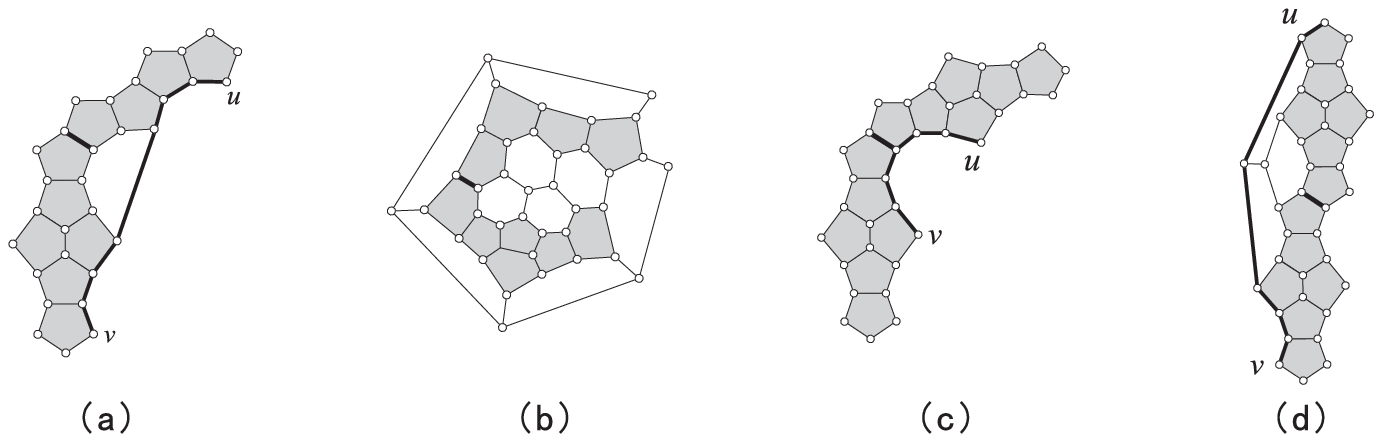}\\
{Figure \ref{fig3-7}: $B_{2}\ast B_{3}$ (grey graphs in (a) and (b))
and $B_{3}^2$ (grey graphs in (c) and (d)).}
\end{center}
\end{figure}

\begin{proof} Up to isomorphism, all cases of $B_{2}\ast B_{3}$ and
$B_{3}^2$ are illustrated as the graphs in grey color in Figure
\ref{fig3-7}. Suppose to the contrary that $B_{2}\ast B_{3},
B_{3}^2\in \mathscr B$. Then the hexagon extension of $B_{2}\ast
B_{3}, B_{3}^2$ are subgraphs of fullerene graphs. So all graphs
illustrated in Figure \ref{fig3-7} are fragments of fullerene
graphs, contradicting either Proposition \ref{pro2-3} or Lemma
\ref{lem2-2}. So $B_{2}\ast B_{3}, B_{3}^2\notin \mathscr B$.
\end{proof}

\begin{figure}[!hbtp]\refstepcounter{figure}\label{fig3-8}
\begin{center}
\includegraphics{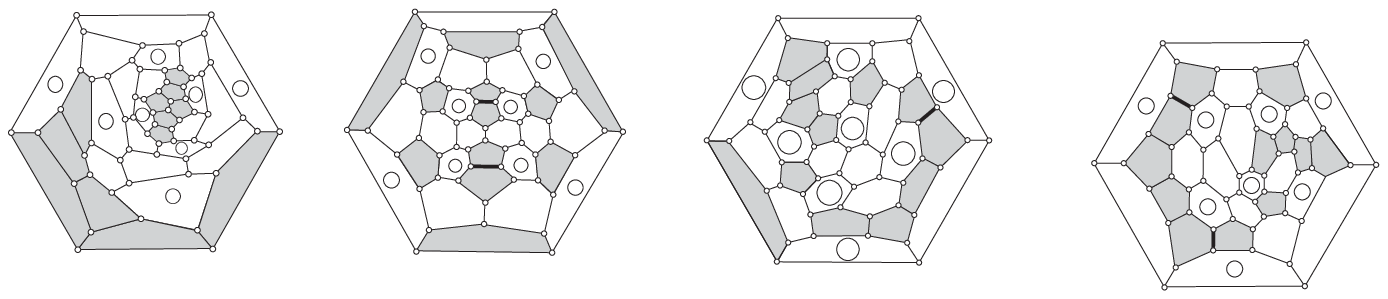}\\
{Figure \ref{fig3-8}: Extremal fullerene graphs with 60 vertices.}
\end{center}
\end{figure}

We are particularly interested in the graphs in $\mathscr B_{60}$.
From the extremal fullerene graphs shown in Figure \ref{fig3-8}, we
can easily see that $\{P, B_{2}, B_{3}, P^2, P\ast B_{2}, P\ast
B_{2}\ast P\}\subseteq \mathscr B_{\ge 60}$. In fact, these two sets
are equal.

\begin{lem}\label{lem3-9}
$\mathscr B_{\ge 60}=\{P, B_{2}, B_{3}, P^2, P\ast B_{2}, P\ast
B_{2}\ast P\}$.
\end{lem}
\begin{proof} It is clear that $\{P, B_{2}, B_{3}, P^2, P\ast B_{2}, P\ast
B_{2}\ast P\}\subseteq \mathscr B_{\ge 60}$. In the following, we
will prove another direction that $\mathscr B_{\ge 60}\subseteq \{P,
B_{2}, B_{3}, P^2, P\ast B_{2}, P\ast B_{2}\ast P\}$. By Lemma
\ref{lem3-8}, it suffices to prove $B_{2}^2\nsubseteq B$ and
$B_{3}\ast P\nsubseteq B$ for any $B\in \mathscr B_{\ge 60}$.

Suppose $B_{2}^2\subseteq B\in \mathscr B_{\ge 60}$. Clearly,
$B_{2}^2$ has two cases as shown in Figure \ref{fig3-9} (the grey
subgraphs in (a) and (c)). Their Clar extensions
$C[B_{2}^2]\subseteq H[B]\subseteq F_n$ are graphs (a) and (c) in
Figure \ref{fig3-9}. The corresponding hexagon extensions
$H[C[B_2^2]]$ are graphs (b) and (d) in Figure \ref{fig3-9}. Since
$H[C[B_2^2]]$ has four 2-degree vertices, $n\le V(T[C[B_2^2]])+2\le
V(H[C[B_2^2]])+2\le 56$ by Lemma \ref{lem2-2} and Proposition
\ref{pro2-4}, contradicting that $n\ge 60$. Hence, $B\notin \mathscr
B_{\ge 60}$.

\begin{figure}[!hbtp]\refstepcounter{figure}\label{fig3-9}
\begin{center}
\includegraphics{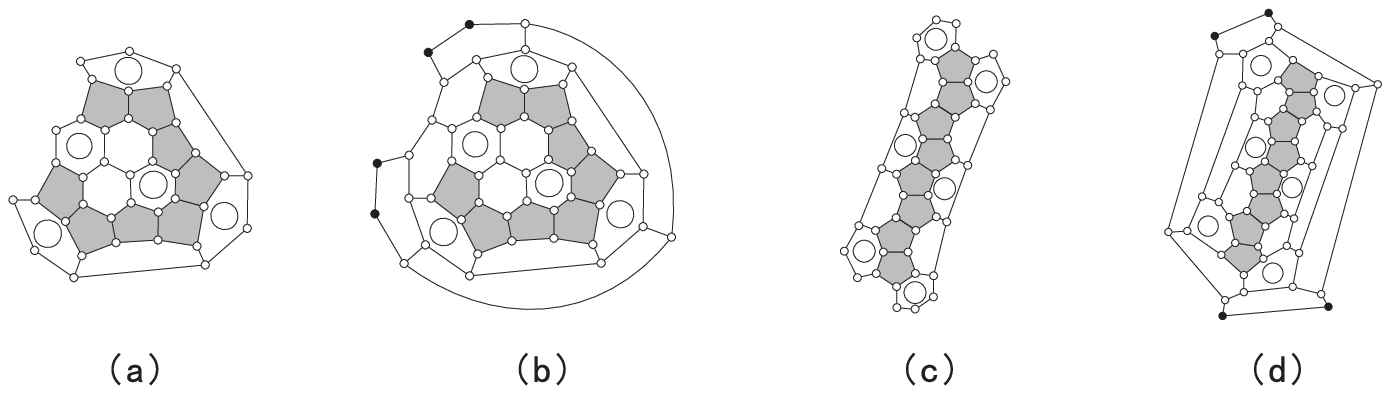}\\
{Figure \ref{fig3-9}: Clar extensions $C[B_{2}^2]$ ((a) and (c)) and
their hexagon extensions ((b) and (d)).}
\end{center}
\end{figure}

Now suppose $B_{3}\ast P\subseteq B \in \mathscr B_{\ge 60}$. Its
Clar extension $C[B_{3}\ast P]\subset H[B]\subseteq F_n$ and
$H[C[B_3\ast P]]$ are illustrated in Figure \ref{fig3-10}. Let $f$
be the face adjoining $H[C[B_3\ast P]]$ as shown in Figure
\ref{fig3-10}. Let $G:=H[C[B_3\ast P]]\cup f$. Then $G$ has at most
four 2-degree vertices. So $n\le |V(G)|+2\le V(H[C[B_3\ast
P]])+3=44$ by Lemma \ref{lem2-2} and Proposition \ref{pro2-4},
contradicting that $n\ge 60$. So $B\notin \mathscr B_{\ge
60}$.\end{proof}

\begin{figure}[!hbtp]\refstepcounter{figure}\label{fig3-10}
\begin{center}
\includegraphics{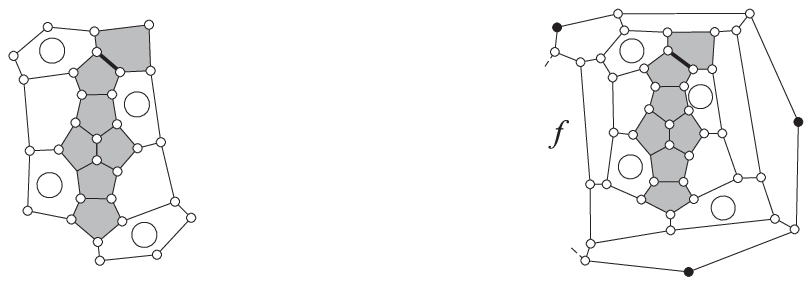}\\
{Figure \ref{fig3-10}: Graphs $C[B_{3}\ast P]$ (left) and
$H[C[B_3\ast P]]$ (right).}
\end{center}
\end{figure}

\begin{thm}\label{thm1}
Let $B$ be a maximal pentagonal fragment of fullerene graph $F_n$
{\upshape ($n\ge 60$)}. Then $B$ is extremal  if and only if $B\in
\mathscr B_{\ge 60}$.
\end{thm}

\begin{proof} By Lemma \ref{lem3-9} and the extremal fullerene graphs in
Figure \ref{fig3-8}, the sufficiency is obvious. So it suffices to
prove the necessity. Let $B$ be a maximal extremal pentagonal
fragment with $k$ pentagons in $F_n$ ($n\ge 60$). We use induction
on $k$ to prove $B\in \mathscr B_{\ge 60}$. Let $\mathcal H$ be a
Clar set of $H[B]$ and $M_B$ be the matching of $B-\mathcal H$ which
covers all remaining 3-degree vertices of $B$. Let $S$ be a subset
of $V(B)$ meeting at most $|S|-1$ pentagons in $B$. If $S\subseteq
U_{\mathcal H}(B)$, then $B-S$ has $k+1-|S|$ pentagons and then has
at least $k+1-|S|$ vertices in $U_{\mathcal H}(B)$ by Lemma
\ref{lem2-7}. Hence $|U_{\mathcal H}(B)|\ge k+1$, contradicting that
$B$ is extremal . So, in the following, we may assume that
$U_{\mathcal H}(B)$ contains no such $S$.

For $k=1$ or $2$, then $B=P$ or $P^2$. The necessity holds since
$P,P^2\in \mathscr B_{\ge 60}$. Now suppose that $k\ge 3$ and the
necessity holds for smaller $k$. Let $p,p_1,p_2$ be the three
pentagons of $B$. By Lemma \ref{lem3-7}, $\gamma (B)=1$. Let $p$ be
the pentagon adjoining only one pentagon, say $p_1$, along an edge
$e$ and $V(p)-V(e)=\{v_1 v_2, v_3\}$. Enumerate clockwise the
hexagons in $F_n$ adjoining $p$ as $h_1, h_2, h_3$ and $h_4$ (see
Figure \ref{fig3-11}). Since any two vertices in $\{v_1,v_2,v_3\}$
form a vertex set $S$, it follows that $U_{\mathcal H}(B)$ contains
at most one of $v_1, v_2$ and $v_3$.

\begin{figure}[!hbtp]\refstepcounter{figure}\label{fig3-11}
\begin{center}
\includegraphics{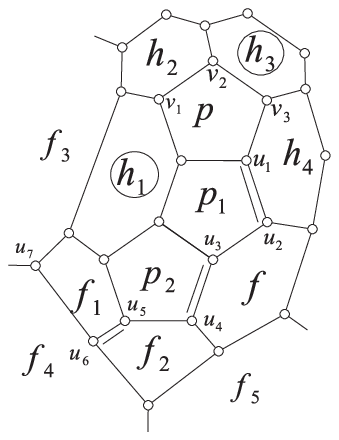}\\
{Figure \ref{fig3-11}: Pentagonal fragment $B$ with $h_1,h_3\in
\mathcal H$.}
\end{center}
\end{figure}

If one of $v_1$ and $v_3$ belongs to $U_{\mathcal H}(B)$, say $v_1$,
then $h_3\in \mathcal H$ since $v_2,v_3\in V(\mathcal H)$. So $h_1,
h_4\notin \mathcal H$. Let $S:=\{v_1\}\cup V(p\cap p_1)$. Then
$S\subseteq U_{\mathcal H}(B)$, contradicting the assumption. If
$v_2\in U_{B}$, then $h_1, h_2\in \mathcal H$. Let
$B':=B-\{v_1,v_2,v_3\}$. Then $B'$ has $k-1$ pentagons and
$|U_{\mathcal H}(B')|=k-1$. By the inductive hypothesis, $B'\in
\mathscr B_{\ge 60}$. So $B$ arises from pasting $B'$ and $p$ along
$p\cap p_1$ and hence $B\in \mathscr B_{\ge 60}$ by Lemma
\ref{lem3-9}. From now on, suppose $v_1,v_2,v_3\in V(\mathcal H)$.

First suppose $h_1,h_3\in \mathcal H$. Let $u_1u_2=p_1\cap h_4$.
Then $u_1u_2\in M_B$. If $f$ is a pentagon, by the symmetry and a
similar discussion as that of Subcase 1.1, we have $B=B_{3}\in
\mathscr B_{\ge 60}$. So suppose $f$ is a hexagon. Then $f\notin
\mathcal H$ since $u_2\in V(f)\cap V(\mathcal H)$. Let
$u_3u_4=p_2\cap f$. Then $u_3u_4\in M_B$. Let $f_1,f_2$ be other two
faces adjoining $p_2$ as illustrated in Figure \ref{fig3-11}. Since
$\{u_1,...,u_4\}\subseteq U_{\mathcal H}(B)$ meets at least four
pentagons, $f_2$ is a pentagon. Let $u_5u_6=f_1\cap f_2$. Then
$u_5u_6\in M_B$ since $f_1\notin \mathcal H$. Let $f_3,f_4,f_5$ be
the other three faces adjoining $f_1$ or $f_2$. Since
$\{u_1,...,u_6\}\subseteq U_{\mathcal H}(B)$ should meet at least
six pentagons, both $f_1$ and $f_4$ are pentagons. Let
$u_6u_7=f_1\cap f_4$. Then $u_7\in U_{\mathcal H}(B)$ because
$f_3\notin \mathcal H$. Since $\{u_1,...,u_7\}\subseteq U_{\mathcal
H}(B)$ meets at least seven pentagons, $f_3$ is also a pentagon. So
$R_5^-=f_1\cup f_2\cup f_3\cup f_4\cup p_2\subseteq B$. By Lemma
\ref{lem3-5}, $B$ is not extremal .

\begin{figure}[!hbtp]\refstepcounter{figure}\label{fig3-12}
\begin{center}
\includegraphics{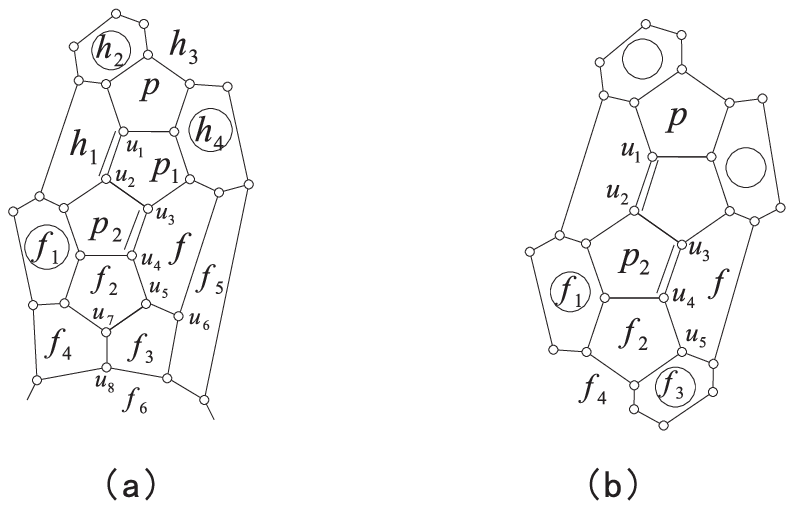}\\
{Figure \ref{fig3-12}: Illustration for the proof of Theorem
\ref{thm1}.}
\end{center}
\end{figure}

So, in the following, suppose $h_2,h_4\in \mathcal H$. Let $f$ be
the face adjoining $p_1, p_2$ and $h_4$, and let $u_1u_2=h_1 \cap
p_1$ and $u_3u_4=f \cap p_2$ (see Figure \ref{fig3-12} (a)). Then
$u_1u_2, u_3u_4\in M_B$.

First suppose $f$ is a hexagon. Let $f_1$ and $f_2$ be the other two
faces adjoining $p_2$, distinct from $h_1, p_1$ and $f$ (see Figure
\ref{fig3-12} (a)). Since $\{u_1, u_2, u_3, u_4\}$ meets at least 4
pentagons, $f_2$ is a pentagon. If $f_1\notin \mathcal H$, then
$S=V(p_2)\subseteq U_{\mathcal H}(B)$, contradicting the assumption.
So $f_1\in \mathcal H$. Let $u_4u_5=f_2\cap f$ and let $f_3,f_4\ne
p_2$ be the other two faces adjoining both $f_2$ as shown in Figure
\ref{fig3-12} (a). If $f_3\notin \mathcal H$, then $f_3$ is a
pentagon since $\{u_1,...,u_5\}\subseteq U_{\mathcal H}(B)$ meets at
least five pentagons. Let $u_5u_6=f\cap f_3$ and $u_7u_8=f_3\cap
f_4$. Clearly, $u_6,u_7\in U_{\mathcal H}(B)$. Let $f_5,f_6$ be two
faces adjoining $f_3$ as shown in Figure \ref{fig3-12} (a). Since
both $\{u_1,u_2,...,u_5,u_6\}$ and $\{u_1,u_2,...,u_5,u_7\}$ meet at
least six pentagons, both $f_5$ and $f_4$ are pentagonal. If $f_6$
is a pentagon, then $f_2,f_3,f_4,f_5$ and $f_6$ form a
$R_5^-\subseteq B$, contradicting Lemma \ref{lem3-5}. So $f_6$ is a
hexagon. Clearly, $f_6\notin \mathcal H$ because $\{u_5u_6,
u_7u_8\}\subset M_B$ or $\{u_5u_7, f_3\cap f_5\} \subset M_B$. So
$S:=V(f_3)\subseteq U_{\mathcal H}(B)$, a contradiction. The
contradiction implies that $f_3\in \mathcal H$. Hence $p\cup p_1\cup
p_2\cup f_2=B_2$ and $f_2\cap f_4$ is a pasting edge (see Figure
\ref{fig3-12} (b)). If $f_4$ is a hexagon, then $B=B_{2}\in \mathscr
B_{\ge 60}$. If $f_4$ is a pentagon, let $B':=B-(p_1\cup p_2\cup
p\cup \{u_5\})$. Then $B'$ has $k-4$ pentagons and $|U_{\mathcal
H}(B')|=k-4$. By inductive hypothesis, $B'\in \mathscr B_{\ge 60}$.
Hence, $B$ arises from pasting $B'$ and $B_{2}$ along $f_2\cap f_4$.
Therefore, $B\in \mathscr B_{\ge 60}$ by Lemma \ref{lem3-9}.

\begin{figure}[!hbtp]\refstepcounter{figure}\label{fig3-13}
\begin{center}
\includegraphics{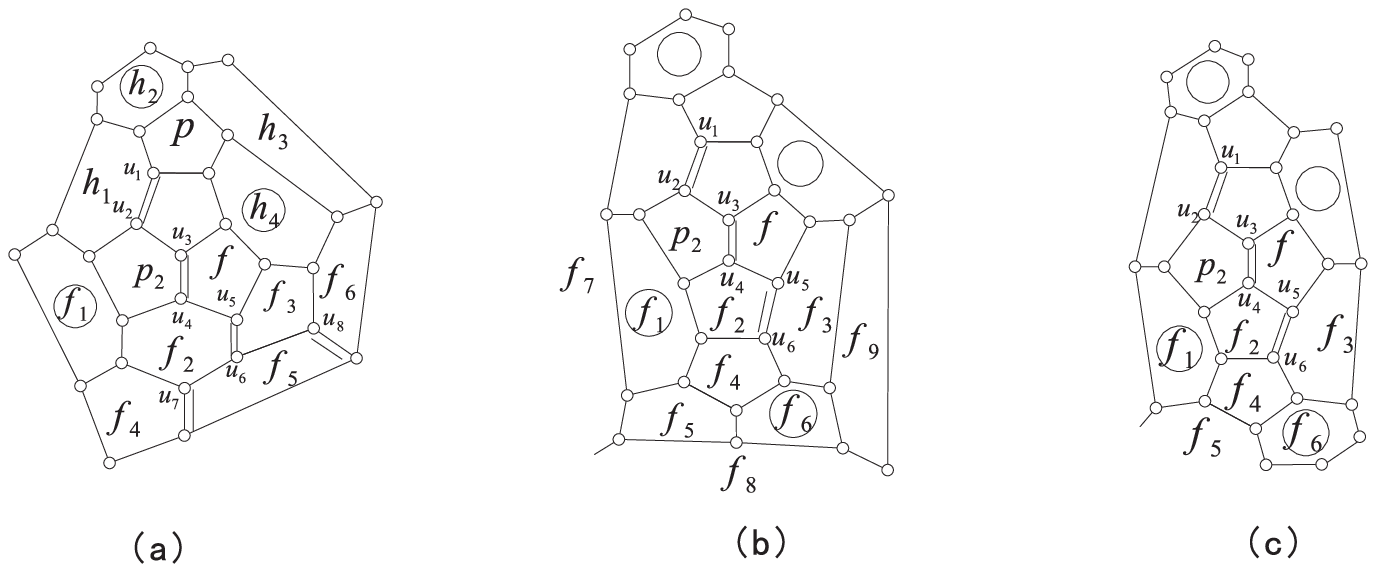}\\
{Figure \ref{fig3-13}: Illustration for the proof of Theorem
\ref{thm1}.}
\end{center}
\end{figure}

Now suppose that $f$ is a pentagon (see Figure \ref{fig3-13}). Let
$u_5u_6=f_2\cap f_3$. Then $f_1$(or $f_3$) and $f_2$ cannot be
pentagonal simultaneously by Lemma \ref{lem3-8}. Since $u_3u_4\in
M_B$, we have $f_2\notin \mathcal H$. Clearly $f_3\notin \mathcal
H$. So $\{u_1,...,u_5\}\subseteq U_{\mathcal H}(B)$. Hence
$\{u_1,...,u_5\}$ meets at least five pentagons. So at least one of
$f_2$ and $f_3$ is pentagonal.

If $f_3$ is a pentagon, then $f_2$ is a hexagon and $u_5u_6\in M_B$
by Lemma \ref{lem3-5}. Let $f_4, f_5$ and $f_6$ be the other three
faces adjoining $f_2$ or $f_3$ as illustrated in Figure
\ref{fig3-13} (a). Let $u_6u_7=f_2\cap f_5$ and $u_6u_8=f_3\cap
f_5$. Since both $\{u_1,...,u_6,u_7\}\subseteq U_{\mathcal H}(B)$
and $\{u_1,...,u_6,u_8\}\subseteq U_{\mathcal H}(B)$ meet at least
seven pentagons, all $f_4, f_5, f_6$ are pentagonal. Hence $f_4\cap
f_5\in M_B$ and $f_5\cap f_6\in M_B$. So $S:=V(f_5)\subseteq
U_{\mathcal H}(B)$ meets only four pentagons in $B$, a
contradiction.

So suppose that $f_2$ is a pentagon and both $f_1$ and $f_3$ are
hexagons (see Figure \ref{fig3-13} (b)). Clearly, $f_3\notin
\mathcal H$ and $u_5u_6\in M_B$ since $h_4\in \mathcal H$. Since
$V(f_2)$ meets four pentagons, $V(f_2)\nsubseteq U_{\mathcal H}(B)$
and hence $f_1\in \mathcal H$. Let $f_4$ be the face adjoining
$f_1,f_2$ and $f_3$. Then $f_4$ is a pentagon since
$\{u_1,...,u_6\}$ meets at least six pentagons. Let $f_5$ and $f_6$
be the faces adjoining $f_4$ as illustrated in Figure \ref{fig3-13}
(b). Clearly, $f_5\notin \mathcal H$ since it is adjacent with
$f_1\in \mathcal H$.

If $f_6\notin \mathcal H$, then $V(f_4\cap f_6)\subset U_{\mathcal
H}(B)$. Both $f_5$ and $f_6$ are pentagons since
$\{u_1,...,u_6\}\cup V(f_4\cap f_6)\subseteq U_{\mathcal H}(B)$
meets at least 8 pentagons. Let $f_7, f_8$ and $f_9$ be faces
adjoining $f_5$ or $f_6$ as illustrated in Figure \ref{fig3-13} (b).
Since $\{u_1,...,u_6\}\cup V(f_3\cap f_6)\subseteq U_{\mathcal
H}(B)$, we have $f_9$ is a pentagon of $B$. Since $\{f_3\cap f_6,
f_5\cap f_6\}\subset M_B$ or $\{f_4\cap f_6, f_6\cap f_9\}\subset
M_B$, we have $f_8\notin \mathcal H$. Further, $f_8$ is a hexagon
because $R_5^-\nsubseteq B$. Hence $S:=V(f_6)\subseteq U_{\mathcal
H}(B)$ meets only four pentagons in $B$, a contradiction.

So suppose that $f_6\in \mathcal H$. Then $p\cup p_1\cup p_2\cup
f\cup f_2\cup f_4=B_{3}$ (see Figure \ref{fig3-13} (c)). By the
proof of Lemma \ref{lem3-8}, we have $B=B_{3}\in \mathscr B_{\ge
60}$. This completes the proof of the theorem. \end{proof}

Let $F_n$ ($n\ge 60$) be an extremal fullerene graph. That means
$c(F_n)=\frac {n-12} 6$. By Lemma \ref{lem3-1} and Theorem
\ref{thm1}, every pentagon of $F_n$ lies in a pentagonal fragment
$B\in \mathscr B_{\ge 60}$. For a Clar formula $\mathcal H$ of $F_n$
and a maximal pentagonal fragment $B$ of $F_n$, we have that
$\mathcal H\cap H[B]$ is a Clar set of $H[B]$ where $H[B]$ is the
hexagon extension of $B$.

\begin{thm}\label{thm2}
Let $F_n$ $(n\ge 60)$ be a fullerene graph and $B_1,B_2,...,B_k$ be
all maximal pentagonal fragments of $F_n$.
Then $F_n$ is extremal if and only if \\
{\upshape (1)} $B_i\in \mathscr B_{\ge 60}$ for all $1\le i\le k$;
and\\
{\upshape (2)} $H[\cup_{i=1}^{k}B_i]$ has a normal Clar set
$\cup_{i=1}^{k}\mathcal H_i$ where $\mathcal H_i$ is the Clar set of
$H[B_i]$; and\\
{\upshape (3)} $F_n-C[\cup_{i=1}^{k}B_i]$ has a sextet pattern
covering all vertices in $V(F_n-C[\cup_{i=1}^{k}B_i])$.
\end{thm}

Theorem \ref{thm2} gives a characterization of extremal fullerne
graphs. This characterization provides an approach to construct all
extremal fullerene graphs with 60 vertices.

\section{Extremal Fullerene graphs with $60$ vertices}

Let $F_n$ be an extremal fullerene graph and $\mathcal H$ be a Clar
formula of $F_{n}$. Then $|\mathcal H|=\frac{n-12} 6$ and
$M:=F_{n}-\mathcal H$ is a matching with six edges. By Theorem
\ref{thm2}, every pentagon lies in a maximal extremal pentagonal
fragment $B \in \mathscr B_{\ge 60}$ and $\mathcal H\cap H[B]$ is a
Clar set of $H[B]$ where $H[B]$ is the hexagon extension of $B$.
Then $M_B=E(B)\cap M$ is the matching of $B$ covering all 3-degree
vertices of $B$ in $V(B-\mathcal H)$. For $B=P^2$, or $B_{2}\ast P$
or $P\ast B_{2}\ast P$, every $P$ has a vertex $v$ uncovered by
$M_B$. Obviously, $v$ is covered by $M$ and let $uv\in M$. Then $u$
belongs to another $P$. The edge $uv$ connects two $P$s to form a
graph $B_{1}$ as illustrated in Figure \ref{fig2-2}.

\begin{pro}\label{pro4-1}
Let $\mathcal H$ be a Clar formula of an extremal fullerene graph
$F_{n}$ $(n\ge 60)$ and $M:=F_{n}-\mathcal H$. Then a face $f$ of
$F_{n}$ is a pentagon if and only if there exists an edge $e\in M$
such that $e\cap f\ne \emptyset$ and $e\notin E(f)$. \qed
\end{pro}

\begin{figure}[!hbtp]\refstepcounter{figure}\label{fig4-1}
\begin{center}
\includegraphics{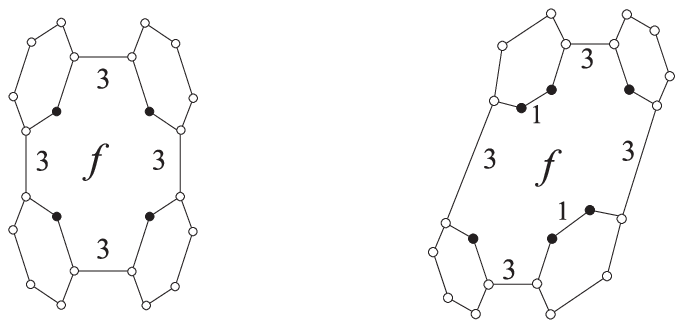}\\
{Figure \ref{fig4-1}: The boundary labelings: $3333$ (right) and
$331331$ (left).}
\end{center}
\end{figure}

Let $G$ be a 2-connected subgraph of $F_n$. Then every face of $G$
is bounded by a cycle. Let $f$ be a face of $G$ with $k$ 2-degree
vertices of $G$. Then $k$ 2-degree vertices separate $f$ into $k$
degree-saturated paths. Use a $k$-length sequence to label $f$ such
that every numbers in the sequence correspond clockwise the lengths
of all degree-saturated paths. The maximum one in the lexicographic
order over all such $k$-length sequences is called the {\em boundary
labeling} of $f$ (see Figure \ref{fig4-1}).

\begin{pro}\label{pro4-2}
Let $B$ be a fragment of an extremal fullerene graph $F_n$ and
$\mathcal H$ be a Clar formula of $F_n$. Let $W$ be the set of all
2-degree vertices on $\partial B$. Then:\\
{\upshape (1)} $|W|\ne 1$;\\
{\upshape (2)} the boundary labeling of $\partial B$ is $ij$ with $5\ge i\ge j\ge 4$ for $|W|=2$;\\
{\upshape (3)} $|W|\ne 3$ for $W\subseteq V(\mathcal H)$;\\
{\upshape (4)} the boundary labeling of $\partial B$ is 3333 or
$i$3$j$1 with $5\ge i\ge j\ge 4$ for $|W|=4$ and $W\subseteq
V(\mathcal H)$.
\end{pro}

\begin{proof} Since $B$ is a fragment, $\partial B$ is a cycle. Let
$C:=\partial B$. For convenience, we may draw $B$ on the plane such
that $C$ bounds an inner face. All 2-degree vertices in $W$ separate
$C$ into $|W|$ degree-saturated paths. Let $v\in W$ and
$vv_1,vv_2\in E(C)$. Let $v_3$ be the third neighbor of $v$ in
$F_n$. Then $v_3$ lies in $F_n-B$ or $W$. Since $F_n$ is
3-connected, $|W|>1$.

If $|W|=2$, then the two 2-degree vertices are adjacent by Lemma
\ref{lem2-2}. Since every face of $F_n$ is either a hexagon or a
pentagon, the length of any degree-saturated path connecting the two
2-degree vertices is either 4 or 5. It follows that the boundary
labeling of $\partial B$ is $ij$ with $5\ge i\ge j\ge 4$.

If $|W|=3$, then the 3-degree vertices have a common neighbor $u$ by
Lemma \ref{lem2-2}. Since $W\subseteq V(\mathcal H)$, it follows
that $u$ is an isolate vertex of $F_n-\mathcal H$, contradicting
that $\mathcal H$ is a Clar formula of $F_n$. So $|W|\ne 3$ if
$W\subseteq V(\mathcal H)$.

\begin{figure}[!hbtp]\refstepcounter{figure}\label{fig4-2}
\begin{center}
\scalebox{0.8}{\includegraphics{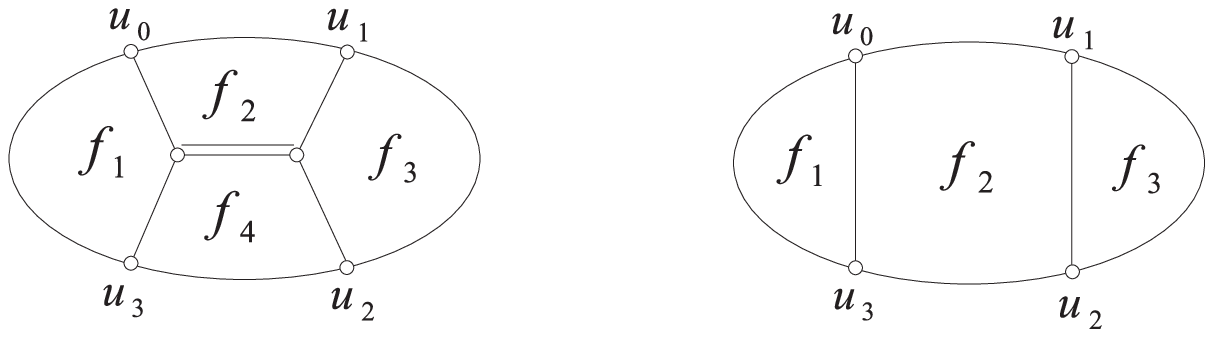}}\\
{Figure \ref{fig4-2}: Illustration for the proof of Proposition
\ref{pro4-2}.}
\end{center}
\end{figure}

Now suppose $|W|=4$. Let $u_0,u_1,u_2,u_3$ be the four vertices
clockwise on $C$ (see Figure \ref{fig4-2}). Let $P_{u_i,u_{i+1}}$
($i, i+1\in \mathbb Z_4$) be the degree-saturated path of $C$
connecting $u_i$ and $u_{i+1}$. Let $T:=F_n-(V(B)\setminus W)$, the
subgraph induced by the vertices within $C$ and the vertices in $W$.
By Lemma \ref{lem2-2}, $T$ is $T_0$ or the union of two $K_2$s or a
3-length path. If $T$ is $T_0$, then the two vertices in the
interior of $C$ are adjacent and hence induce an edge $e$. Then
$e\in M:=F_n-\mathcal H$. Let $f_1,f_2,f_3,f_4$ be the four faces
meeting the edge $e$ (see Figure \ref{fig4-2} (left)). By
Proposition \ref{pro4-1}, both $f_1$ and $f_3$ are pentagonal. So
$|P_{u_3,u_0}|=4$ and $|P_{u_1,u_2}|=4$. Since $5\le |f_3|\le 6$ and
$5\le |f_4|\le 6$, we have that $u_0u_1\notin E(F_n)$ and
$u_2u_3\notin E(F_n)$. Then $u_0$ and $u_1$ cannot be in the common
hexagon in $\mathcal H$. Similarly, $u_2$ and $u_3$ cannot be in the
common hexagon in $\mathcal H$. So $|P_{u_0u_1}|=4$ and
$|P_{u_2u_3}|=4$. Hence all $P_{u_i,u_{i+1}}$ for $i, i+1\in \mathbb
Z_4$ are 3-length path. Further, the boundary labeling of $\partial
B$ is 3333.

If $G$ is the union of two $K_2$s or a 3-length path, then $u_0u_3,
u_1u_2\in E(F_n)$. Let $f_1,f_2,f_3$ be the three faces of $F_n$
within $C$ (see Figure \ref{fig4-2} (right)). Hence $5\le
|P_{u_0u_1}|\le 6$ and $5\le |P_{u_2,u_3}|\le 6$ since $5\le
|f_1|\le 6$ and $5\le |f_3|\le 6$. Since $5\le |f_2|\le 6$ and
$\{u_0,u_1,u_2,u_3\}\subseteq V(\mathcal H)$, then one of
$|P_{u_0,u_1}|$ and $|P_{u_2,u_3}|$ equals 2 and the other equals 4.
It follows that the boundary labeling of $\partial B$ is $i3j1$ with
$5\ge i\ge j\ge 4$. \end{proof}

In the following, $F_{60}$ always means an extremal fullerene graph
with 60 vertices. Using $B_1$ instead of $P$ in the pasting
operation, let $\mathscr G_{60}$ denote the set of all maximal
subgraphs of $F_{60}$ arising from the pasting operation on $B_1,
B_2$ and $B_3$. Up to isomorphism, the Clar extension of $G\in
\mathscr G_{60}$ is unique since the Clar extension of any element
in $\mathscr B_{\ge 60}$ is unique. Note that $B_1^k$ is the graph
obtained by pasting $k$ graphs isomorphic to $B_1$ along the pasting
edge of each $P$ in $B_1$.

\begin{lem}\label{lem4-3}
$\mathscr G_{60}\subseteq \{B_1, B_{2}, B_{3}, B_1^2, B_1^3, B_1^4,
B_{1}\ast B_{2}, B_{1}\ast B_{2}\ast B_1, B_1\ast B_{2}\ast
B_{1}\ast B_{2}, B_{2}\ast B_{1}\ast B_{2}\}.$
\end{lem}

\begin{proof} Since $c(F_{60})=8$, we have that $B_1^k$ and $(B_1\ast
B_{2})^r$ satisfy $k\le 4$ and $r\le 2$ if they belong to $\mathscr
G_{60}$.

By Lemma \ref{lem3-9}, it suffices to prove $B_1^2\ast
B_{2}\nsubseteq G$ for any $G\in \mathscr G_{60}$. Suppose to the
contrary that $B_1^2\ast B_{2}\subseteq G\in \mathscr G_{60}$. Then
either $G=B_1^3\ast B_{2}$ or $G=B_1^2\ast B_{2}$ by $c(F_{60})=8$.

\begin{figure}[!hbtp]\refstepcounter{figure}\label{fig4-3}
\begin{center}
\includegraphics{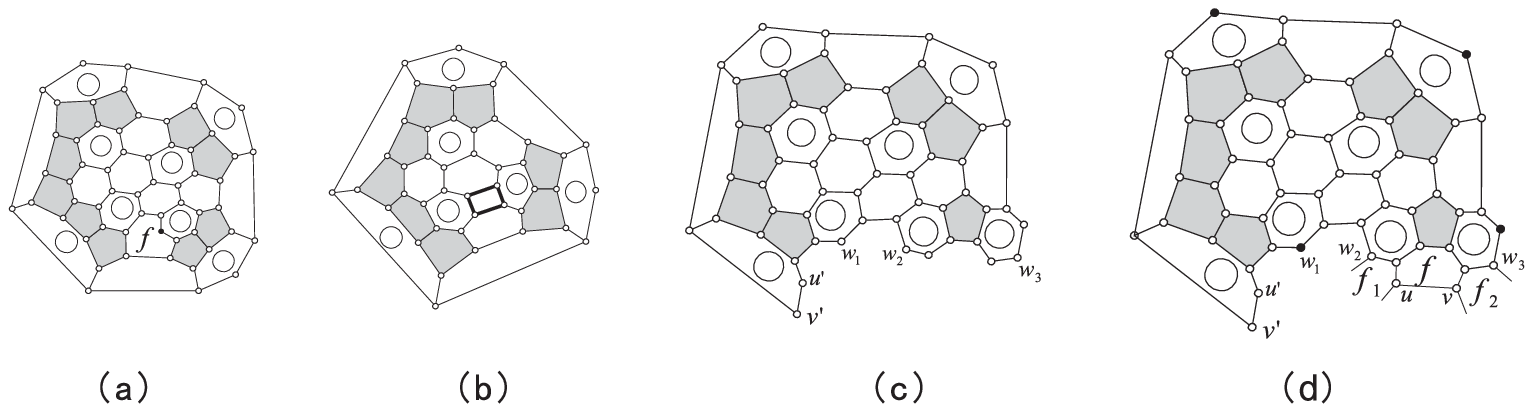}\\
{Figure \ref{fig4-3}: Illustration for the proof of Lemma
\ref{lem4-3}.}
\end{center}
\end{figure}

If $G=B_1^3\ast B_{2}$, then $B_1^3\ast B_{2}$ has to be the grey
subgraph of the graph (a) in Figure \ref{fig4-3} since
$c(F_{60})=8$. The subgraph induced by $C[G]$ in $F_{60}$ is the
graph (a) in Figure \ref{fig4-3}. Proposition \ref{pro4-2} implies
that the graph (a) is not a subgraph of $F_{60}$. Hence $B_1^3\ast
B_{2}\notin \mathscr G_{60}$, a contradiction.

If $G=B_1^2\ast B_{2}$, then there are two cases for $G$ as the grey
subgraphs illustrated in graphs (b) and (c) in Figure \ref{fig4-3},
respectively. The graphs (b) and (c) are the subgraphs induced by
$C[G]$. Clearly, the graph (b) could not be a subgraph of $F_n$ in
that it has a 4-length cycle. For the graph (c), let $f$ be the
hexagon adjoining $G$ along an edge of $B_1$ and $u,v,u',v',w_1,
w_2, w_3$ be some 2-degree vertices on the boundary of $G\cup f$
(see Figure \ref{fig4-3} (d)). If $uv\in E(\mathcal H)$, then $u=u'$
and $v=v'$ since $F_{60}$ is a cubic plane graph. Then $w_1$ is
adjacent to $w_2$ by Lemma \ref{lem2-2}, which forms a 4-length
cycle in $F_{60}$, a contradiction. So suppose $uv\in M$. Let $f_1$
and $f_2$ be the pentagons met by $uv$ but not containing it by
Proposition \ref{pro4-1}. Whether $uv\in M_{B_1}$ or $uv\in
M_{B_2}$, one of $f_1$ and $f_2$ adjoins two hexagons in $\mathcal
H$. So either $u'v'\in E(f_1)$ or $u'v'\in E(f_2)$. If $u'v'\in
E(f_1)$, then $w_2$ is adjacent to $u'$ and hence $w_1$ would be a
unique 2-degree on a face of a subgraph of $F_{60}$, contradicting
Proposition \ref{pro4-2}. So suppose $u'v'\in E(f_2)$. Then $w_3$ is
adjacent to $v'$, which forms a face with three 2-degree vertices
which belong to $V(\mathcal H\cap C[G])$, also contradicting
Proposition \ref{pro4-2}. So $B_1^2\ast B_{2}\notin \mathscr
G_{60}$. This completes the proof. \end{proof}

\begin{figure}[!hbtp]\refstepcounter{figure}\label{fig4-4}
\begin{center}
\scalebox{0.9}{\includegraphics{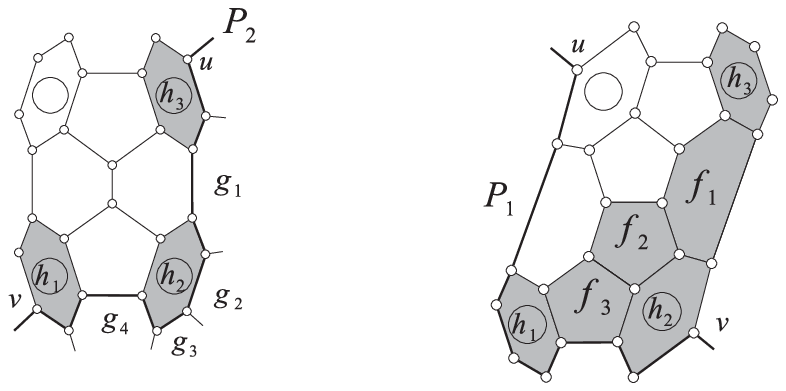}}\\
{Figure \ref{fig4-4}: Clar extensions of $B_1$ and $B_2$.}
\end{center}
\end{figure}

\begin{lem}\label{lem4-4}
Let $G\subset F_{60}$ such that $G$ has two components, one of which
is $B_1$ and another is $B_1$ or $B_{2}$. If the Clar extension
$C[G]$ of $G$ is a fragment, then $|C[G]\cap \mathcal H|\ge 6$.
\end{lem}

\begin{proof} Let $B_1$ and $B$ be two components of $G$, where $B$
is isomorphic to $B_1$ or $B_2$. By Theorem \ref{thm2}, the Clar set
of $C[G]$ is a subset of a Clar formula $\mathcal H$ of $F_{60}$.
Clearly, $|C[B_1]\cap \mathcal H|=4$ and $|C[B]\cap \mathcal H|=4$.
Then $|C[G]\cap \mathcal H|=|(C[B_1]\cap \mathcal H)\cup (C[B]\cap
\mathcal H)|-|C[B_1]\cap C[B]\cap \mathcal H|$. If $|C[B_1]\cap
C[B]\cap \mathcal H|\le 2$, then $|C[G]\cap \mathcal H|\ge 6$ and
the lemma is true.

So suppose $|C[B_1]\cap C[B]\cap \mathcal H|\ge 3$ and let $h_1,
h_2,h_3\in C[B_1]\cap C[B]\cap \mathcal H$ (see Figure
\ref{fig4-4}). Let $B'\subset C[B]$ be a fragment such that $B'$
contains $h_1, h_2, h_3$ and has minimal number of inner faces. Then
$B'$ has at most 6 inner faces including $h_1,h_2$ and $h_3$ (see
Figure \ref{fig4-4}, the faces $f_1,f_2,f_3$ in $C[B_{2}]$) and
$B'\cap C[B_1]=h_1\cup h_2\cup h_3$. Since $C[G]$ is a fragment, the
faces of $B'$ different from $h_1,h_2,h_3$ adjoins $C[B_1]$. It
needs at least 4 faces adjoining $C[B_1]$ to join $h_1, h_2$ and
$h_3$ to form a fragment (the faces $g_1,...,g_4$ in $C[B_1]$, see
Figure \ref{fig4-4}). So $B'$ has at least 7 inner faces,
contradicting that $B'$ has at most 6 faces. The contradiction
implies that $|C[B_1]\cap C[B]\cap \mathcal H|\le 2$. So the lemma
is true. \end{proof}

\begin{lem}\label{lem4-5}
If $B_{3}\subset F_{60}$, then $F_{60}$ contains no other elements
in $\mathscr{G}_{60}$ as subgraphs.
\end{lem}
\begin{proof} Let $H[B_{3}]$ be the hexagon extension of $B_{3}$. Then
$H[B_{3}]\subset F_{60}$. Let $f_1$ and $f_2$ be the two hexagons
adjoining $B_{3}$ and let $f_3,f_4$ be two faces adjoining $C[B_3]$
as shown in Figure \ref{fig4-5} (a).

\begin{figure}[!hbtp]\refstepcounter{figure}\label{fig4-5}
\begin{center}
\includegraphics{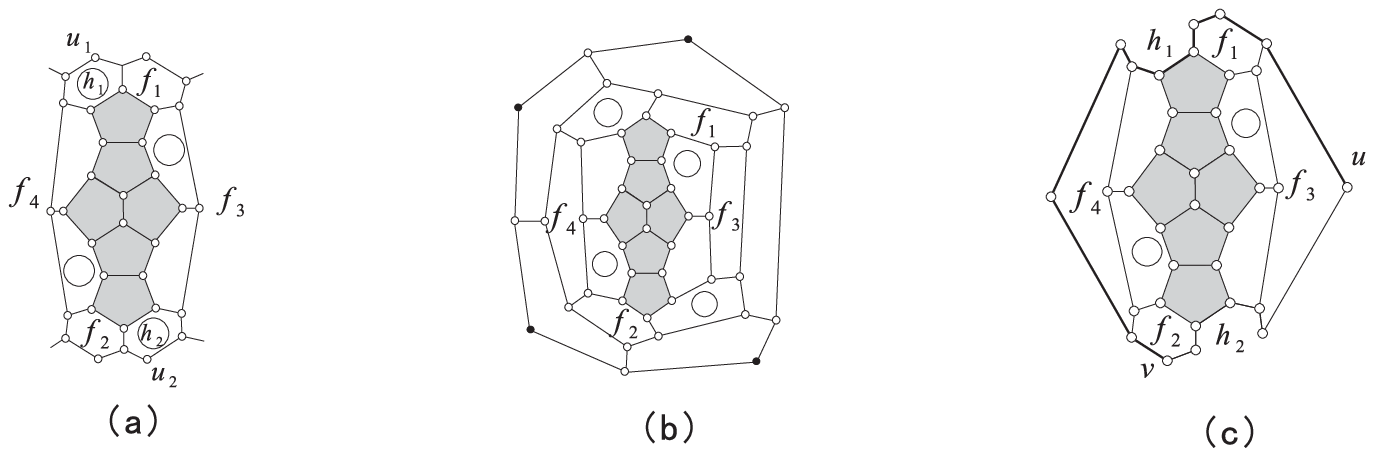}\\
{Figure \ref{fig4-5}: Illustration for the proof of Lemma
\ref{lem4-5}.}
\end{center}
\end{figure}

Let $G_1:=H[B_{3}]\cup f_3\cup f_4$. If at least one of $f_3$ and
$f_4$, say $f_3$, is a pentagon. Then the hexagon extension $H[G_1]$
of $G_1$ contains at most four 2-degree vertices on its boundary
(see Figure \ref{fig4-5} (b)). By Lemma \ref{lem2-2} and Proposition
\ref{pro2-4}, it holds that $n\le |V(G_1)|+9+2\le 46$ if $G_1\subset
F_n$. So suppose both $f_3$ and $f_4$ are hexagons since $G_1\subset
F_{60}$.

Let $h_1, h_2\in \mathcal H\cap H[B_{3}]$ and  $u_i\in V(h_i)$
($i=1,2$) as illustrated in Figure \ref{fig4-5} (a). Let
$G_2:=G_1-\{u_1, u_3\}$ (see Figure \ref{fig4-5} (c)). By
Proposition \ref{pro4-1}, we have $V(\partial G_2)\subset V(\mathcal
H)$. Let $G_3:=F_{60}-(G_2-\partial G_2)$. Then $G_3$ has six
pentagons and $|G_3\cap\mathcal H|=6$. Let $f$ be the unique face of
$G_3$ which is not a face of $F_{60}$. Then $G_3\cup G_2=F_{60}$ and
$G_3\cap G_2=f=\partial G_2$. So a 2-degree vertex (resp. 3-degree
vertex) of $G$ on $f$ is identified to a 3-degree vertex (resp.
2-degree vertex) on $\partial G_2$ in $F_{60}$.

If $B_{3}\nsubseteq G_3$, then every hexagon in $G_3\cap \mathcal H$
belongs to either $C[B_1]$ or $C[B_{2}]$. For $h_i\in G_3\cap
\mathcal H$ ($i=1,2$), let $P_i=\partial C[B]\cap \partial G_2$
where $B=B_1$ or $B_{2}$. Since $F_{60}$ is cubic, $|P_i|\ge 11$ for
$i=1,2$ (the thick paths on $\partial C[B_1]$ or $\partial C[B_{2}]$
connecting vertices $u$ and $v$ in Figure \ref{fig4-4}). Therefore,
$|V(f)|\ge 11+11-2=20$ which contradicts $|V(f)|=|V(\partial
G_2)|=16$. So $B_{3}\subset G_3$. \end{proof}

\begin{lem}\label{lem4-6}
There are two distinct extremal fullerene graphs which have 60
vertices and contain $B_{3}$ as subgraphs.
\end{lem}
\begin{proof} If $B_{3}\subset F_{60}$, then $F_{60}$ contains two subgraphs
isomorphic to $B_3$ by Lemma \ref{lem4-5}. Let $C[B_{3}]$ be the
Clar extension of $B_{3}$ (see Figure \ref{fig4-6} (a)). If two
subgraphs isomorphic to $C[B_{3}]$ have common hexagons in $\mathcal
H$, according to the proof of Lemma \ref{lem4-5}, the common
hexagons belong to $\{h_1, h_2\}$ (see Figure \ref{fig4-6} (a)). By
the symmetry, let $h_2$ be a common hexagon. Let $f_1,f_2$ be two
faces adjoining the $C[B_{3}]$ as shown in Figure \ref{fig4-6} (b).
Then one of $f_1$ and $f_2$ is a pentagon of the second $B_{3}$
since $h_2$ belongs to the Clar set of the second $C[B_{3}]$. If
$f_1$ is a pentagon, then a fullerene graph $F_{48}$ is formed as
illustrated in Figure \ref{fig4-6} (b). If $f_2$ is a pentagon, then
another fullerene graph $F_{48}$ is formed as illustrated in Figure
\ref{fig4-6} (c).

\begin{figure}[!hbtp]\refstepcounter{figure}\label{fig4-6}
\begin{center}
\scalebox{0.9}{\includegraphics{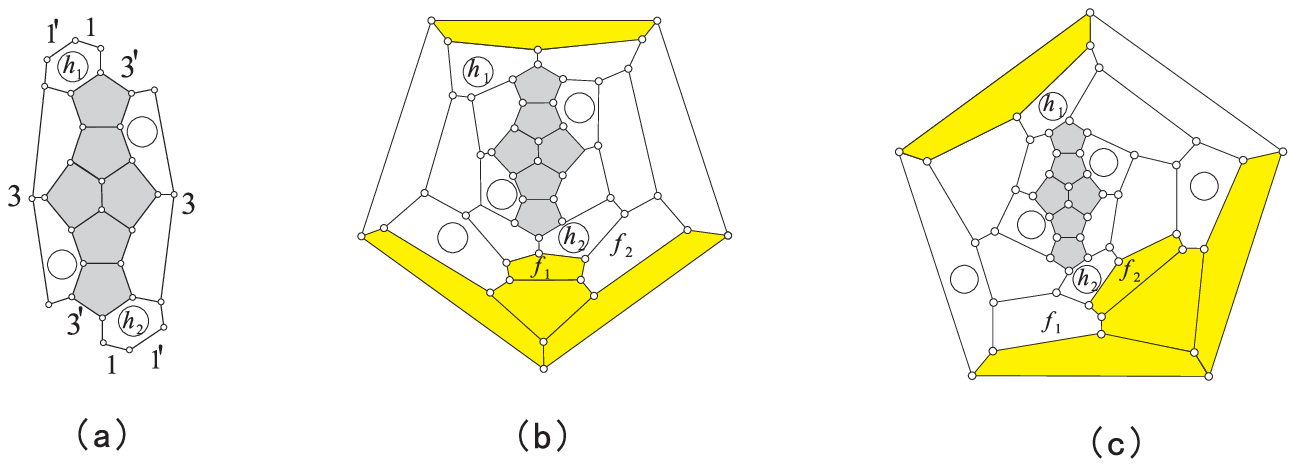}}\\
{Figure \ref{fig4-6}: Illustration for the proof of Lemma
\ref{lem4-6}.}
\end{center}
\end{figure}

So suppose the two subgraphs isomorphic to $C[B_{3}]$ have no common
hexagon in $\mathcal H$. Further, the two subgraphs isomorphic to
$C[B_{3}]$ have no common vertex. Since $|V(C[B_{3}])|=30$, hence
$F_{60}$ is formed by using edges to connect the 2-degree vertices
on the boundaries of the two subgraphs isomorphic to $C[B_{3}]$. On
the other hand, the faces of $F_{60}$ do not belong to two
$C[B_{3}]$s are hexagons. The boundary labeling of $C[B_{3}]$ is
33113311. Hence the 3-length degree-saturated path of one $C[B_{3}]$
together with the 1-length degree-saturated path of another
$C[B_{3}]$ form a hexagon. Since the paths with same length have two
distinct positions on the $\partial C[B_{3}]$, use the labeling
$33'11'33'11'$ (see Figure \ref{fig4-6} (a)) to distinguish the same
length degree-saturated paths with different positions. If the new
hexagons consist of either the paths with label 3 and the paths with
label 1 or the paths with label $3'$ and the paths with label $1'$,
then a $F_{60}$ is formed as illustrated in Figure \ref{fig4-7}
(left). If the new hexagons consists of either the paths with label
$3'$ and the paths with label 1 or the paths with label 3 and the
paths with label $1'$, then another $F_{60}$ is formed as
illustrated in Figure \ref{fig4-7} (right). So there are exactly two
extremal fullerene graphs $F_{60}^1$ and $F_{60}^2$ with $B_{3}$ as
subgraphs. \end{proof}

\begin{figure}[!hbtp]\refstepcounter{figure}\label{fig4-7}
\begin{center}
\includegraphics{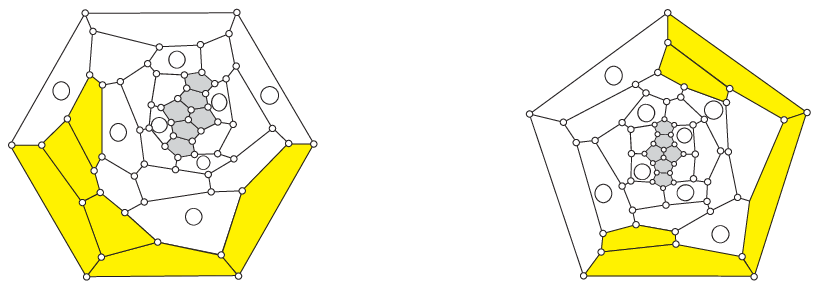}\\
{Figure \ref{fig4-7}: Extremal fullerene graphs $F^1_{60}$ and
$F_{60}^2$.}
\end{center}
\end{figure}

\begin{lem}\label{lem4-7}
There are six distinct extremal fullerene graphs which have 60
vertices and contain $B_1^k$ $(2\leq k\leq 4)$ as subgraphs.
\end{lem}
\begin{proof} \noindent{\em Case 1:} $B_1^4\subset F_{60}$ is maximal. Since
$c(F_{60})=8$, we have that $B_1^4$ is unique and its Clar extension
$C[B_1^4]$ is the graph illustrated in Figure \ref{fig4-8} (a). By
Lemma \ref{lem2-2}, we have two different extermal fullerene graphs
$F^3_{60}$ and $F^4_{60}$ as shown in Figure \ref{fig4-8} (b) and
(c).

\begin{figure}[!hbtp]\refstepcounter{figure}\label{fig4-8}
\begin{center}
\includegraphics{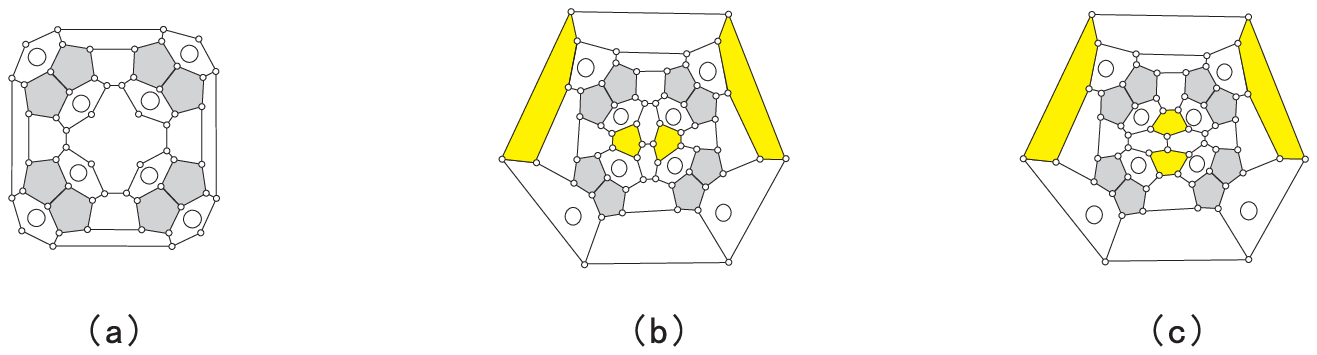}\\
{Figure \ref{fig4-8}: Extremal fullerene graphs $F^3_{60}$ and
$F^4_{60}$ with $B_1^4$ as maximal subgraphs.}
\end{center}
\end{figure}

\noindent{\em Case 2:} $B_{1}^3\subset F_{60}$ is maximal. There are
two cases for $B_1^3$ whose Clar extensions are illustrated in
Figure \ref{fig4-9} (a) and (b). By Proposition \ref{pro4-2}, the
graph (a) is not a subgraph of $F_{60}$. So $B_1^3\subset F_{60}$ is
unique and its Clar extension $C[B_1^3]$ is the graph (b). Let
$h_1,h_2,...,h_8$ be the all eight hexagons in $C[B_1^3]\cap
\mathcal H$ and let $v,v_1,v_2,...,v_7,u,u_1,u_2,...,u_7$ be all
2-degree vertices on the boundary of $C[B_1^2]$ as shown in Figure
\ref{fig4-9} (b). Let $f_1$ be the face adjoining $C[B_1^2]$ (see
Figure \ref{fig4-9} (b)).

\begin{figure}[!hbtp]\refstepcounter{figure}\label{fig4-9}
\begin{center}
\includegraphics{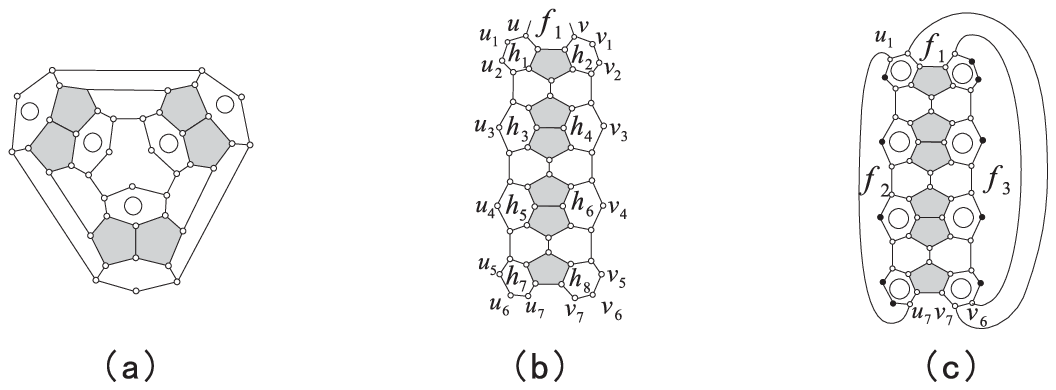}\\
{Figure \ref{fig4-9}: Illustration for the proof of Case 2.}
\end{center}
\end{figure}

If $f_1$ adjoins three hexagons in $\mathcal H$, then either
$v_5v_6\in E(f_1)$ or $v_6v_7\in E(f_1)$ by symmetry. If $v_5v_6\in
E(f_1)$, then $v,v_1$ are adjacent to $v_5,v_4$, respectively. Then,
by Lemma \ref{lem2-2}, $v_2$ is adjacent to $v_3$. Then the edge
$v_2v_3$ together with the 3-length degree-saturated path connecting
$v_2$ and $v_3$ form a 4-length cycle in $F_{60}$, a contradiction.
So suppose $v_6v_7\in E(f_1)$. Then $u_1$ is adjacent to $u_7$ (see
Figure \ref{fig4-9} (c)). Let $f_2$ and $f_3$ be the two faces
adjoining the graph (c). Each of $f_2$ and $f_3$ has five 2-degree
vertices. Let $I[f_i]$ ($i=2,3$) be the subgraph consisting of $f_i$
together with its interior. Then $I[f_2]$ and $I[f_3]$ together
contain three edges in $M$. One of them, say $I[f_2]$, satisfies
that $I[f_2]-f_2$ is an edge in $M$. However, the two ends of one
edge are adjacent to at most four 2-degree vertices on $f_2$ since
$F_{60}$ is cubic, contradicting that $f_2$ has five 2-degree
vertices.

\begin{figure}[!hbtp]\refstepcounter{figure}\label{fig4-10}
\begin{center}
\includegraphics{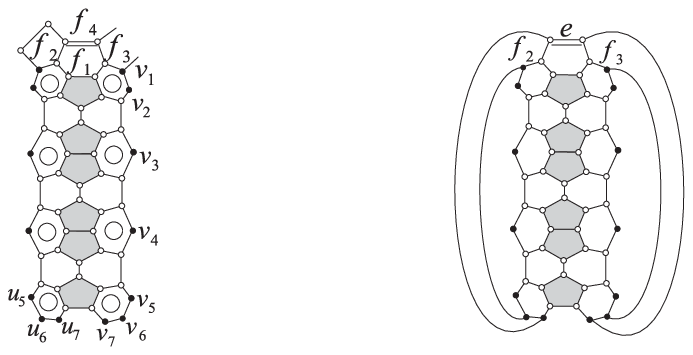}\\
{Figure \ref{fig4-10}: Illustration for the proof of Case 2.}
\end{center}
\end{figure}

So suppose that $f_1$ contains an edge $e=w_1w_2\in M$. By
proposition \ref{pro4-1}, let $f_2$ and $f_3$ be two pentagons such
that $f_2\cap f_1=w_1u$ and $f_3\cap f_1=w_2v$ (see Figure
\ref{fig4-10} (left)). According to Lemma \ref{lem4-5}, either $e\in
E(B_{2})\cap M$ or $e\in E(B_1)\cap M$. First suppose $e\in
E(B_{2})$. Let $f_4$ be the pentagon containing $e$ (see Figure
\ref{fig4-10} (left)). Then one of $f_2$ and $f_3$, say $f_3$, is
adjacent to two hexagons in $\mathcal H$. Then $f_3$ adjoins either
$h_8$ or $h_7$ since $c(F_{60})=8$. Note that $v_6v_7\notin E(f_2)$
and $u_5u_6\notin E(f_3)$ since $f_4$ is a pentagon. So suppose
either $v_5v_6\in E(f_2)$ or $u_6u_7\in E(f_2)$. If $v_5v_6\in
E(f_3)$, then $v_1$ is adjacent to $v_5$ and hence $v_2,v_3,v_4\in
V(\mathcal H)$ are the all 2-degree vertices on a face boundary,
contradicting Proposition \ref{pro4-2}. If $u_6u_7\in E(f_3)$, then
$v_1$ is adjacent to $u_7$ and hence $v_2, v_3$ are adjacent to
$v_7,v_6$, respectively. Furthermore, $v_4$ is adjacent to $v_5$ by
Lemma \ref{lem2-2}. Hence a subgraph of $F_{60}$ with a 4-length
cycle is formed, a contradiction.

So suppose $e\in E(B_1)$. Then both $f_2$ and $f_3$ adjoin two
hexagons in $\mathcal H$. Hence, $f_2$ and $f_3$ adjoin $h_7$ and
$h_8$, respectively. Obviously, $v_6v_7\in E(f_2)$ and $u_6u_7\in
E(f_3)$ (see Figure \ref{fig4-10} (right)). By Lemma \ref{lem2-2},
there are three distinct extremal fullerene graphs $F^5_{60}$, and
$F^6_{60}$ and $F^7_{60}$ with the graph as shown in Figure
\ref{fig4-10} (right) as a subgraph (see Figure \ref{fig4-11}).

\begin{figure}[!hbtp]\refstepcounter{figure}\label{fig4-11}
\begin{center}
\includegraphics{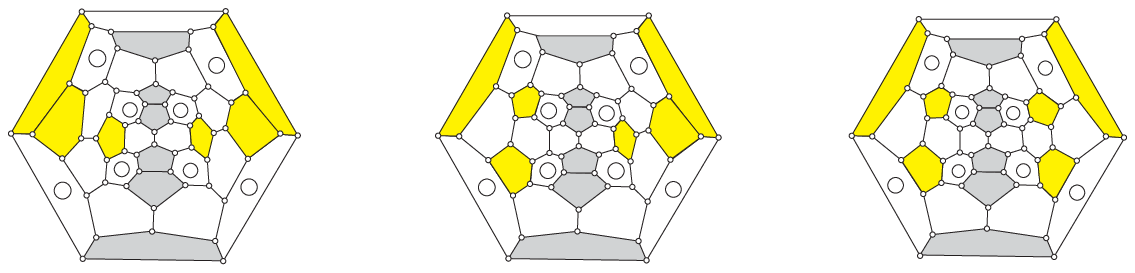}\\
{Figure \ref{fig4-11}: Extremal fullerene graphs  $F^5_{60}$,
$F^6_{60}$ and $F^7_{60}$.}
\end{center}
\end{figure}

\noindent{\em Case 3:} $B_1^2\subset F_{60}$ is maximal and
$B_1^3\nsubseteq F_{60}$. Then $|C[B_1]\cap \mathcal H|\ge 6$. Let
$h_1,...,h_6$ be the six hexagons in $C[B_1]\cap \mathcal H$ as
illustrated in Figure \ref{fig4-12} (a). Let $v_1,...,v_7$ and
$u_1,...,u_7$ be the all 2-degree vertices on the $\partial
C[B_1^2]$ and let $f_1,f_2$ be two hexagons adjoining $C[B_1^2]$
such that $u_1,v_1\in V(f_1)$ and $u_7,v_7\in V(f_2)$ (see Figure
\ref{fig4-12} (a)). Obviously, $f_1\ne f_2$. Let $uv\in E(f_1)$,
then either $uv\in M$ or $uv\in E(\mathcal H)$.

\begin{figure}[!hbtp]\refstepcounter{figure}\label{fig4-12}
\begin{center}
\includegraphics{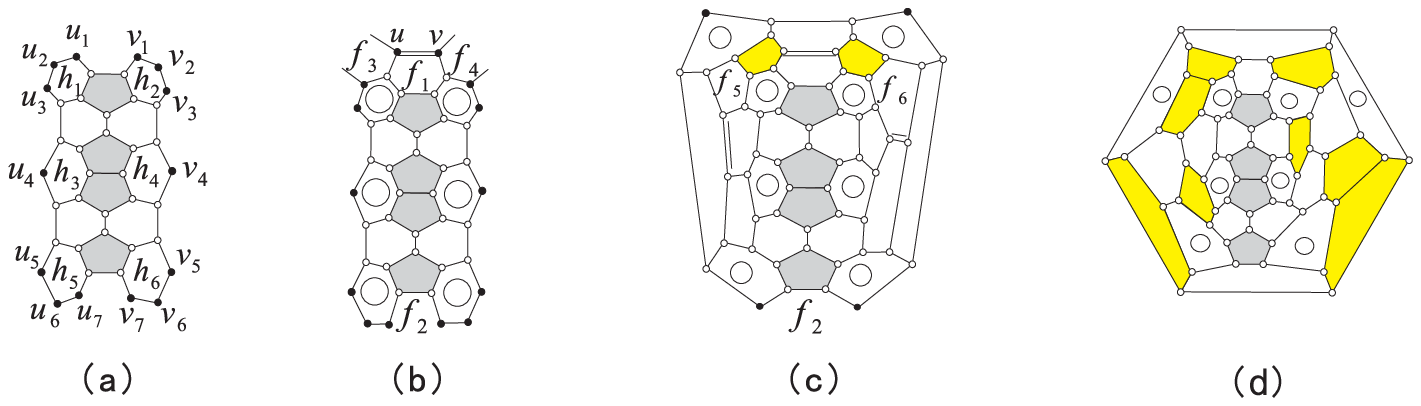}\\
{Figure \ref{fig4-12}: Illustration for the proof of Case 3 and the
extremal fullerene graph $F^8_{60}$.}
\end{center}
\end{figure}

If $uv\in M$, then either $uv\in E(B_1)$ or $uv\in E(B_{2})$. By
Proposition \ref{pro4-1}, let $f_3$ and $f_4$ be the pentagons met
by $uv$ but not containing it (see Figure \ref{fig4-12} (b)). If
$uv\in E(B_1)$, by Proposition \ref{pro4-2}, $f_3$ does not adjoin
$h_5$. If $f_3$ adjoins $h_6$, then either $v_5v_6\in E(f_3)$ or
$v_6v_7\in E(f_3)$. If $v_5v_6\in E(f_3)$, then $v$ is adjacent to
$v_5$ and hence $v$ is adjacent to $v_4$ to bound a hexagon. Then a
subgraph of $F_{60}$ is formed, which has a face with only $v_2,v_3$
connected by a 1-length degree-saturated path on its boundary,
contradicting Proposition \ref{pro4-2}. So suppose $v_6v_7\in
E(f_3)$. Then $u_2$ is adjacent to $v_7$ and hence $u_3$ is adjacent
to $u_7$. A subgraph of $F_{60}$ is formed, which has a face with
only three 2-degree vertices $u_4,u_5,u_6\in V(\mathcal H)$, also
contradicting Proposition \ref{pro4-2}. By symmetry, $f_4$ does not
adjoin $h_5$ and $h_6$. So $f_3$ and $f_4$ adjoin the two hexagons
in $\mathcal H\setminus \{h_1,...,h_6\}$ (see Figure \ref{fig4-12}
(c)).

\begin{figure}[!hbtp]\refstepcounter{figure}\label{fig4-13}
\begin{center}
\includegraphics{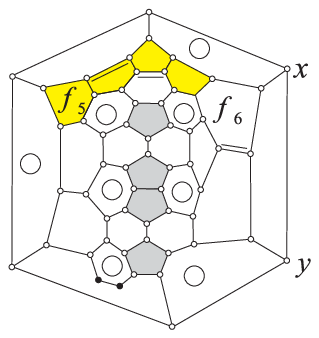}\\
{Figure \ref{fig4-13}: Illustration for the proof of Case 3.}
\end{center}
\end{figure}

Let $f_5$ and $f_6$ be the faces adjoining the $B_1$ with $uv\in
E(B_1)$ along its pasting edges (see Figure \ref{fig4-12}). Since
$B_1^3\nsubseteq F_{60}$, at least one of $f_5$ and $f_6$ is a
hexagon. If both $f_5$ and $f_6$ are pentagonal, then $F_{60}$ still
contains a $B_1^3$ which contains three edges in $M$ as $M\cap
E(f_5), M\cap E(f_6)$ and $M\cap E(f_2)$ since $f_2$ is a hexagon.
By symmetry, we may assume $f_5$ is a pentagon and $f_6$ is a
hexagon. Then we have a graph as illustrated in Figure \ref{fig4-12}
(c) which has four 2-degree vertices on its boundary. By Lemma
\ref{lem2-2}, there is a unique extremal fullerene graph $F^8_{60}$
which contains three subgraphs isomorphic to $B_1^2$ as maximal
subgraphs (see Figure \ref{fig4-12} (d)) since $f_2$ is hexagon. Now
suppose $uv\in E(B_{2})$. Let $f_5$ and $f_6$ be the faces adjoining
$f_3$ and $f_4$, respectively. By symmetry, say $f_5\subset B_{2}$
(see Figure \ref{fig4-13}). By Proposition \ref{pro4-2}, the Clar
extension of the $B_{2}$ containing $f_5$ has two hexagons in
$\mathcal H\setminus \{h_1,...,h_6\}$. Whether $f_6$ is a hexagon or
a pentagon, the vertex $x$ is adjacent to $y$ in $F_{60}$. Hence a
subgraph of $F_{60}$ is formed as the graph in Figure \ref{fig4-13},
contradicting Proposition \ref{pro4-2}.

\begin{figure}[!hbtp]\refstepcounter{figure}\label{fig4-14}
\begin{center}
\includegraphics{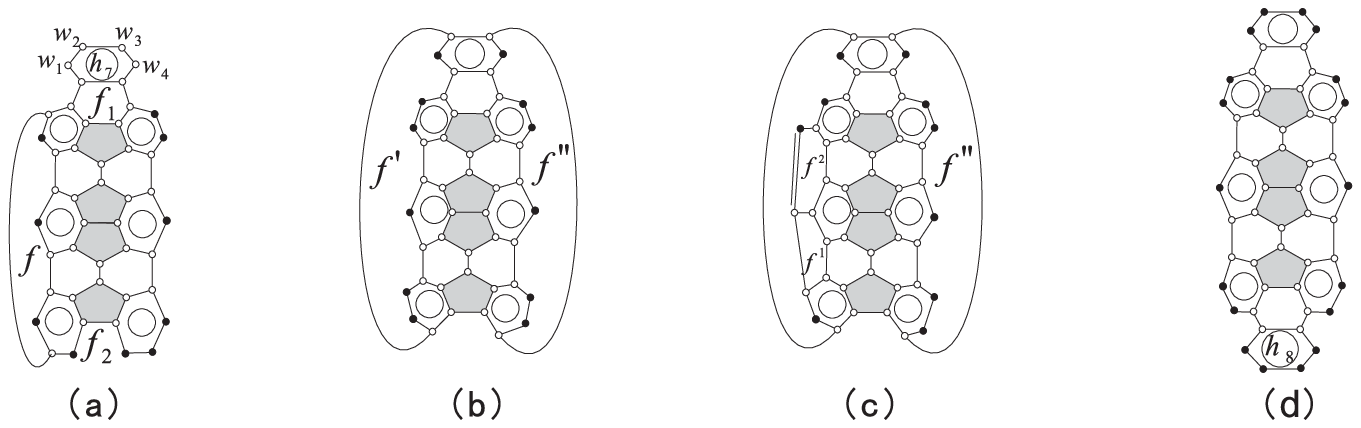}\\
{Figure \ref{fig4-14}: Illustration for the proof of Case 3.}
\end{center}
\end{figure}

So, in the following, suppose that $uv\in E(\mathcal H)$. By
Proposition \ref{pro4-2}, $uv\notin E(h_5)$ and $uv\notin E(h_6)$.
Let $h_7\in \mathcal H$ and $uv\in E(h_7)$ and let the vertices of
$h_7-uv$ be $w_1,w_2,w_3$ and $w_4$ (see Figure \ref{fig4-14} (a)).
By the symmetry of $f_1$ and $f_2$, assume $f_2$ also adjoins three
hexagons in $\mathcal H$. Then $f_2$ adjoins only hexagons $\mathcal
H\setminus\{h_1,h_2,...,h_6\}$. If $f_2$ adjoins $h_7$, then either
$w_1w_2\in E(f_2)$ or $w_2w_3\in E(f_2)$ by symmetry of $w_1w_2$ and
$w_3w_4$. If $w_1w_2\in E(f_2)$, then $w_1$ and $u_2$ are adjacent
to $u_7$ and $u_6$, respectively. Therefore, a subgraph of $F_{60}$
with a face $f$ with three 2-degree vertices $u_3,u_4,u_5\in
V(\mathcal H)$ is formed (see Figure \ref{fig4-14} (a)),
contradicting Proposition \ref{pro4-2}. So suppose $w_2w_3\in
E(f_2)$, then $w_2$ and $w_3$ are adjacent to $u_7$ and $v_7$,
respectively. Let $f'$ and $f''$ be two faces as illustrated in the
graph (b) in Figure \ref{fig4-14}. By symmetry of $f'$ and $f''$, we
may assume that the unique hexagon $\mathcal H\setminus
\{h_1,h_2,...,h_7\}$ lies in the $f'$. Let $I[f']$ and $O[f']$ be
the subgraphs of $F_{60}$ consisting of $f'$ together with it
interior and $f'$ together with its exterior, respectively. Let
$f^1,f^2,f^3,f^4$ be the four faces of $F_{60}$ adjoining $O[f']$
along the four 3-length degree-saturated paths. If one of them is a
pentagon, say $f^1$, then $f^1$ contains a vertex covered by one
edge $e\in M$. Let $e\in E(f_2)$ (see Figure \ref{fig4-14} (c)).
Then $O[f']\cup f^1\cup f^2$ has a face with only four 2-degree
vertices on its boundary. Note that the hexagon in $\mathcal H\cap
I[f']$ has to lie within this face, contradicting that $c\lambda
(F_{60})=5$. So all face of $f^i$ ($i=1,2,3,4$) are hexagons. Let
$G\in \mathscr G_{60}$ lie within $I[f']$. Then $G$ contains or
adjoins at most three hexagons in $\mathcal H\setminus
\{h_2,h_4,h_6,h_3,h_7\}$, which contradicts that a Clar set of
$H[G']$ has at least four hexagons for any $G'\in \mathscr G_{60}$.

So suppose $f_2$ adjoins the hexagon $h_8\in \mathcal H\setminus
\{h_1,h_2,...,h_7\}$ (see Figure \ref{fig4-14} (d)). Let $G''$ be
the graph (d) in Figure \ref{fig4-14}. Then $G''$ contains all
hexagons in $\mathcal H$. So all eight vertices of $F_{60}-V(G'')$
are covered by four edges in $M$ which belong to $E(B_1)$ or
$E_(B_2)$. That means joining some 2-degree vertices on $\partial
G''$ will forming some faces with boundary labeling 3333
(corresponding to the inner face of $C[B_1]-M$ with 2-degree
vertices) or 331331 (corresponding to the inner face of $C[B_2]-M$
with 2-degree vertices) (see Figure \ref{fig4-1}). Hence the
boundary labeling of $\partial G''$ should contains $13331$ or
$1313311$ as subsequences, which contradicts the boundary labeling
of $\partial G''$ is $331311131331311131$. So $G''\nsubseteq
F_{60}$.

Combining Cases 1, 2 and 3, we have exact six fullerene graphs
$F_{60}$ which contain $B_I^k$ $(2\leq k\leq 4)$ as subgraphs.
\end{proof}

\begin{lem}\label{lem4-8}
There are four distinct extremal fullerene graphs $F_{60}$ such that
$B_{2}\ast B_1\subset F_{60}$ and $B_1^k\nsubseteq F_{60}$ $(2\le
k\le 4)$.
\end{lem}
\begin{proof} \noindent{\em Case 1:} $B_{2}*B_{1}*B_{2}*B_{1}\subset F_{60}$
is maximal. Then $B_{2}*B_{1}*B_{2}*B_{1}$ has two different cases
as illustrated in Figure \ref{fig4-15} (a) and (c) since
$c(F_{60})=8$. The Clar extension of the graph (a) induces an
extremal fullerene graph $F^9_{60}$ as shown in Figure \ref{fig4-15}
(b). By Proposition \ref{pro2-3}, the graph (c) is not a subgraph of
$F_{60}$. So there exists a unique $F_{60}$ containing
$B_{2}*B_{1}*B_{2}*B_{1}$ as a subgraph.

\begin{figure}[!hbtp]\refstepcounter{figure}\label{fig4-15}
\begin{center}
\includegraphics{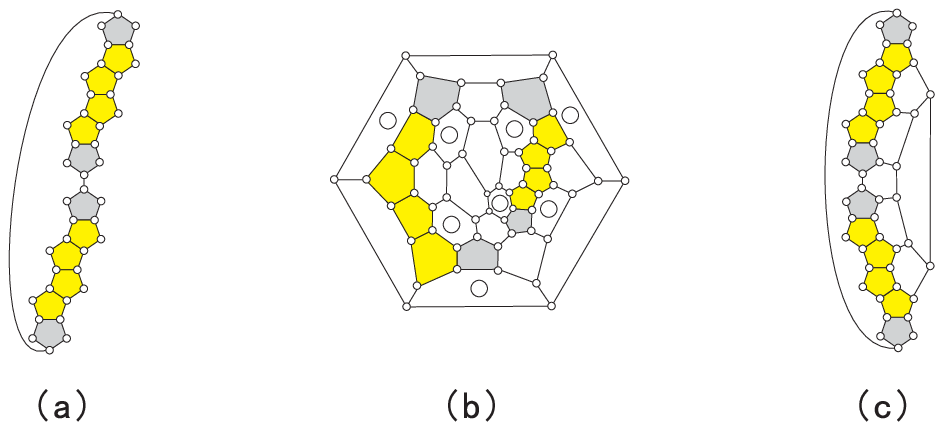}\\
{Figure \ref{fig4-15}: Illustration for the proof of Case 1 and the
extremal fullerene graph $F^9_{60}$.}
\end{center}
\end{figure}

\noindent{\em Case 2:} $B_{2}*B_{1}*B_{2}\subset F_{60}$ is maximal.
By Lemma \ref{lem4-3}, $B_{2}^2\nsubseteq F_{60}$. So
$B_{2}*B_{1}*B_{2}$ has two different cases as illustrated in Figure
\ref{fig4-16} (a) and (c). Their Clar extension induces the graphs
(b) and (d). Both the graphs (b) and (d) have a face $f$ with four
2-degree vertices on its boundary. By Lemma \ref{lem2-2}, an
extremal fullerene graph containing the graph (b) has Clar number
seven. So the graph (b) is not a subgraph of $F_{60}$. From the
graph (d), only one fullerene graph $F^{10}_{60}$ contains
$B_2*B_1*B_2$ as a maximal subgraph (see Figure \ref{fig4-16} (e)).

\begin{figure}[!hbtp]\refstepcounter{figure}\label{fig4-16}
\begin{center}
\includegraphics{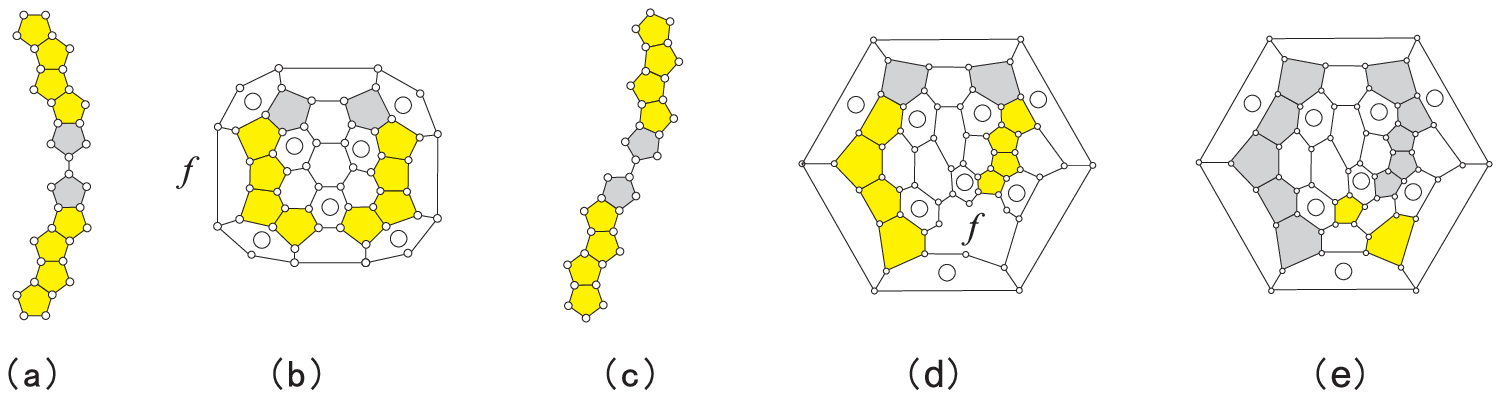}\\
{Figure \ref{fig4-16}: Illustration for the proof of Case 2 and the
extremal fullerene graph $F^{10}_{60}$.}
\end{center}
\end{figure}

\noindent{\em Case 3:} $B_{1}*B_{2}*B_{1}\subset F_{60}$ is maximal.
By the proof of Lemma \ref{lem4-3}, $B_{1}*B_{2}*B_{1}$ is unique as
shown in Figure \ref{fig4-17} (a). Its Clar extension induces the
graph (b), which has a face $f$ with six 2-degree vertices on its
boundary. So the remaining four pentagons adjoin at most four
hexagons in $\mathcal H$ which are adjacent with $f$ in the graph
(b). Hence the four pentagons belong to a $B_{2}$ by Lemma
\ref{lem4-4}. So there is a unique fullerene graph $F^{11}_{60}$
contains $B_1*B_2*B_1$ as a maximal subgraph (see Figure
\ref{fig4-17} (c)).

\begin{figure}[!hbtp]\refstepcounter{figure}\label{fig4-17}
\begin{center}
\includegraphics{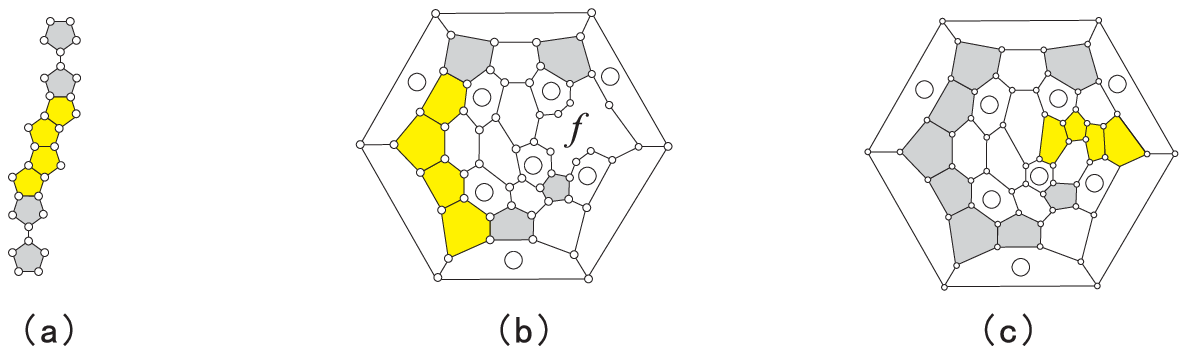}\\
{Figure \ref{fig4-17}: Illustration for the proof of Case 3and the
extremal fullerene graph $F^{11}_{60}$.}
\end{center}
\end{figure}

\noindent{\em Case 4:} $B_1\ast B_{2}\subset F_{60}$ is maximal.
Then it is unique as shown in Figure \ref{fig4-18} (a). Let
$f_1,f_2$ be two faces adjoining the Clar extension $C[B_1\ast B_2]$
as shown in Figure \ref{fig4-18} (a). Since $B_1\ast B_{2}$ is
maximal, $f_1$ and $f_2$ are two pentagons. By Proposition
\ref{pro4-2}, $f_1$ contains an edge $e$ such that $e\notin
E(C[B_1\ast B_{2}])$ and $e\notin E(f_2)$. Clearly, $e\in M$ or
$e\in E(\mathcal H)$.

\begin{figure}[!hbtp]\refstepcounter{figure}\label{fig4-18}
\begin{center}
\scalebox{0.95}{\includegraphics{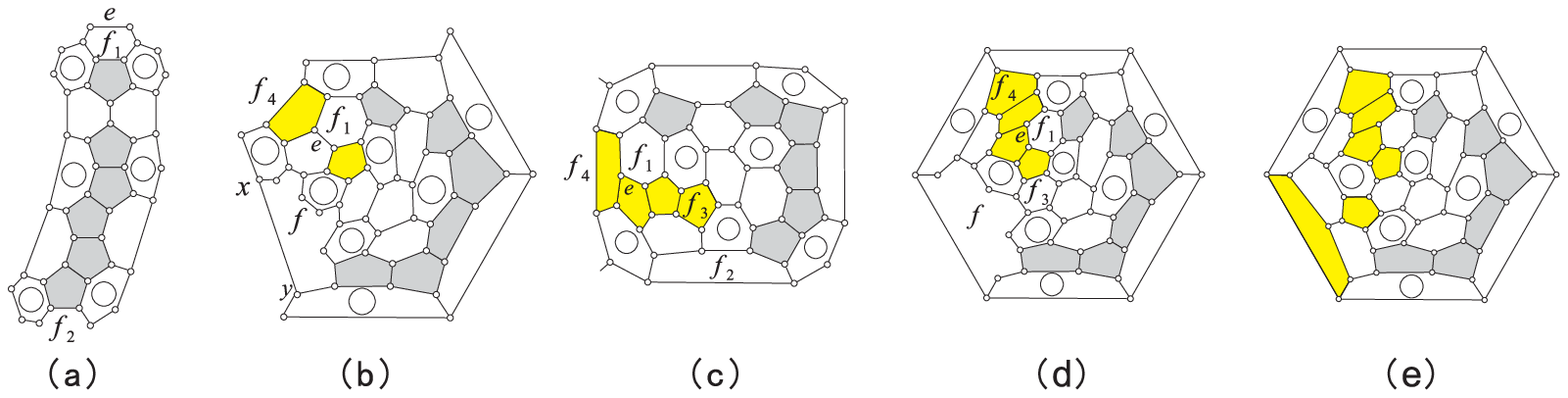}}\\
{Figure \ref{fig4-18}: Illustration for the proof of Subcase 4.1 and
the extremal fullerene graph $F^{12}_{60}$.}
\end{center}
\end{figure}

\noindent{\em Subcase 4.1:} $e\in M$. By Lemma \ref{lem4-5}, either
$e\in E(B_{2})$ or $e\in E(B_1)$.

If $e\in E(B_1)$, then the $B_1$ adjoins two new hexagons in
$\mathcal H$ by Proposition \ref{pro4-2}. Let $x, y\in V(\mathcal
H)$ as shown in Figure \ref{fig4-18} (b). Then the $C[B_1]\cup
C[B_1\ast B_{2}]$ is the graph (e) without the edge $xy$ in Figure
\ref{fig4-18}. Whether $f_4$ is a pentagon or a hexagon, $x$ is
always adjacent to $y$. Hence, a subgraph of $F_{60}$ is formed,
which has a face $f$ with four 2-degree vertices in $V(\mathcal H)$
and with boundary labeling 5313, contradicting Proposition
\ref{pro4-2}.

So suppose $e\in E(B_{2})$. All faces meeting $e$ except $f_1$ are
pentagonal. Let $f_3$ and $f_4$ be the faces adjoining $C[B_1\ast
B_{2}]$ as shown in Figure \ref{fig4-18} (c) and (d). Then either
$f_3$ is a pentagon of $B_{2}$ or $f_4$ is a pentagon of $B_{2}$. If
$f_3$ is a pentagon, then the $C[B_2]\cup C[B_1\ast B_2]$ is the
graph (c) in Figure \ref{fig4-18}. Since the $C[B_2]\cup C[B_1\ast
B_2]$ has four 2-degree vertices on its boundary and has only seven
hexagons in $\mathcal H$, it is not a subgraph $F_{60}$ by Lemma
\ref{lem2-2}. So suppose $f_4$ is a pentagon of the $B_2$. Then the
$C[B_{2}]\cup C[B_1\ast B_2]$ is the graph (d) in Figure
\ref{fig4-18}. By Lemma \ref{lem2-2} and that $B_1\ast B_2$ is
maximal in $F_{60}$, there is a unique fullerene graph $F^{12}_{60}$
containing the graph (d) (see figure \ref{fig4-18} (e)).

\noindent{\em Subcase 4.2:} $e\in E(\mathcal H)$. Let $h\in \mathcal
H$ be the hexagon such that $e\in E(h)$. By Proposition
\ref{pro4-2}, $f_2$ contains an edge $e'$ such that $e'\notin
E(C[B_1\ast B_{2}]\cup h)$ (see Figure \ref{fig4-19} (a)). Then
either $e'\in M$ or $e'\in E(\mathcal H)$.

\begin{figure}[!hbtp]\refstepcounter{figure}\label{fig4-19}
\begin{center}
\includegraphics{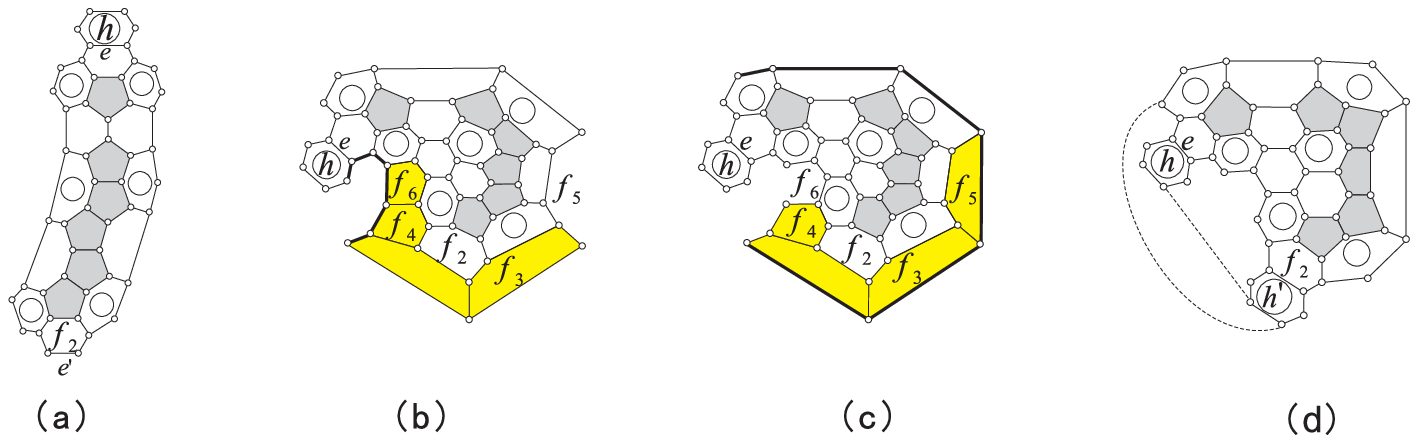}\\
{Figure \ref{fig4-19}: Illustration for the proof of Subcase 4.2.}
\end{center}
\end{figure}

If $e'\in M$, then either $e'\in E(B_1)$ or $e'\in E(B_{2})$ by
Lemma \ref{lem4-5}. If $e'\in E(B_1)$, by Proposition \ref{pro4-2},
then the Clar extension $C[B_1]$ contains two hexagons in $\mathcal
H$ which are different from the seven hexagons in $\mathcal H\cap
(C[B_1\ast B_{2}]\cup h)$. Further, $|\mathcal H|\ge 9$ contradicts
$c(F_{60})=8$. So suppose $e'\in E(B_{2})$. Let $f_3,f_4$ be the two
pentagons meeting $e'$ but $e'\notin E(f_3\cup f_4)$. Let $f_5$ and
$f_6$ be two faces adjoining $f_3$ and $f_4$, respectively (see
Figure \ref{fig4-19} (b) and (c)). Whether $f_5\subset B_{2}$ or
$f_6\subset B_{2}$, we always have a fragment with a 6-length
degree-saturated path on its boundary (see Figure \ref{fig4-19} (b)
and (c), the thick paths), contradicting Proposition \ref{pro2-3}.

So suppose $e'\in E(\mathcal H)$. Let $h'\in \mathcal H$ be the
hexagon containing $e'$ and different from the seven hexagons in
$C[B_1\ast B_2]\cup h$. Let $G$ be the graph induced by $C[B_1\ast
B_2]\cup h\cup h'$ (the graph (d) in Figure \ref{fig4-19}, without
broken lines). Its boundary labeling is 33313111333111 and all
2-degree vertices on it belong to $V(\mathcal H)$. If $G\subset
F_{60}$, then the six vertices in $V(F_{60})\setminus V(G)$ are
covered by three edges in $M\setminus (M\cap E(G))$ and belong to a
$B_1$ or a $B_2$ by Lemma \ref{lem4-5}. So joining some 2-degree
vertices on the boundary of the graph (d) will from some faces with
boundary labeling 3333 (corresponding to $C[B_1]-M$) or 331331
(corresponding to $C[B_{2}]$). That means that the boundary labeling
of $\partial G$ should contain $13331$ (corresponding to $C[B_1]-M$)
or $1313311$ (corresponding to $C[B_2]-M$) as subsequences. Clearly,
33313111333111 contains two subsequences $13331$. So joining four
2-degree vertices on $\partial G$ by two edges will form two faces
with boundary labeling 3333 (see Figure \ref{fig4-19} (d), the dash
edges). Hence, we have a subgraph of $F_{60}$ with a face
(containing the two dash edges) which has a 7-length
degree-saturated path, contradicting Proposition \ref{pro2-3}. So
there is no $F_{60}$ containing $G$.

Combing Cases 1, 2, 3 and 4, there are four extremal fullerene
graphs $F_{60}$ which contain $B_{2}\ast B_1$ as a maximal subgraph
and do not contain $B_1^k$ for $2\le k\le 4$. \end{proof}

\begin{lem}\label{lem4-9}
There are six distinct fullerene graphs $F_{60}$ such that any
$B_1\subset F_{60}$ and any $B_{2}\subset F_{60}$ are maximal.
\end{lem}

\begin{proof} It is well known that $\text{C}_{60}$ is the unique fullerene
graph with 60 vertices and without adjoining pentagons. So
$\text{C}_{60}$ is the unique $F_{60}$ with six subgraphs isomorphic
to $B_1$ as maximal subgraphs (see Figure \ref{fig1-1}). So if
$F_{60}\ne \text{C}_{60}$, then $B_{2}\subset F_{60}$. Let $f_1$ and
$f_2$ be the  two hexagons in the hexagon extension $H[B_{2}]$ and
let $e_i\in E(f_i)$ ($i=1,2$) (see Figure \ref{fig4-20} (a)). It is
easy to see $f_1\cap f_2=\emptyset$ and hence $e_1\ne e_2$. Then
either $e_i\in M$ or $e_i\in E(\mathcal H)$.

\begin{figure}[!hbtp]\refstepcounter{figure}\label{fig4-20}
\begin{center}
\scalebox{0.9}{\includegraphics{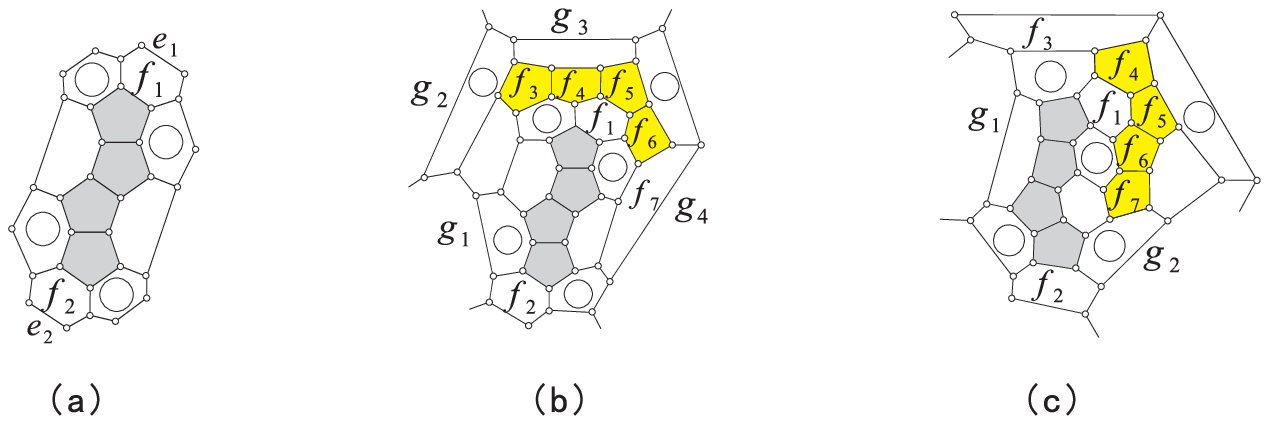}}\\
{Figure \ref{fig4-20}: Illustration for the proof of Case 1.}
\end{center}
\end{figure}

\noindent{\em Case 1:} $e_1,e_2\in M$. Let $f_3, f_4, f_5, f_6$ and
$f_7$ be the faces adjoining $H[B_2]$ as shown in Figure
\ref{fig4-20} (b) and (c). If $e_1$ belongs to a subgraph isomorphic
to $B_2$, denote it by $B'_2$ to distinguish it from the $B_2$ in
Figure \ref{fig4-20} (a). Then either $B'_2=\cup_{i=3}^6 f_i$ or
$B'_2=\cup_{i=4}^7 f_i$. If the former holds, then the $C[B_{2}]\cup
C[B'_2]$ induces the graph (b) in Figure \ref{fig4-20}. Let
$g_1,g_2,g_3,g_4$ be the faces adjoining $C[B_{2}]\cup C[B'_2]$ as
illustrated in Figure \ref{fig4-20}. Note that $C[B_{2}]\cup
C[B'_2]\cup g_1\cup g_3\cup g_4$ has at most four 2-degree vertices
on its boundary. By Lemma \ref{lem2-2}, if $C[B_{2}]\cup
C[B'_2]\subset F_n$, then $n\le 52$. So suppose $B'_2=\cup_{i=4}^7
f_i$. Then $f_3$ is hexagon since $B'_{2}$ is maximal. The
$C[B_{2}]\cup C[B'_2]$ induces the graph (c) in Figure
\ref{fig4-20}. Let $g_1, g_2$ adjoin $C[B_{2}]\cup C[B'_2]$ as shown
in Figure \ref{fig4-20} (c). Then $C[B_{2}]\cup C[B'_2]\cup g_1\cup
g_2$ has at most four 2-degree vertices on its boundary. By Lemma
\ref{lem2-2}, we have $n\le 46$ if $C[B_{2}]\cup C[B'_2]\subset
F_n$.

\begin{figure}[!hbtp]\refstepcounter{figure}\label{fig4-21}
\begin{center}
\scalebox{0.9}{\includegraphics{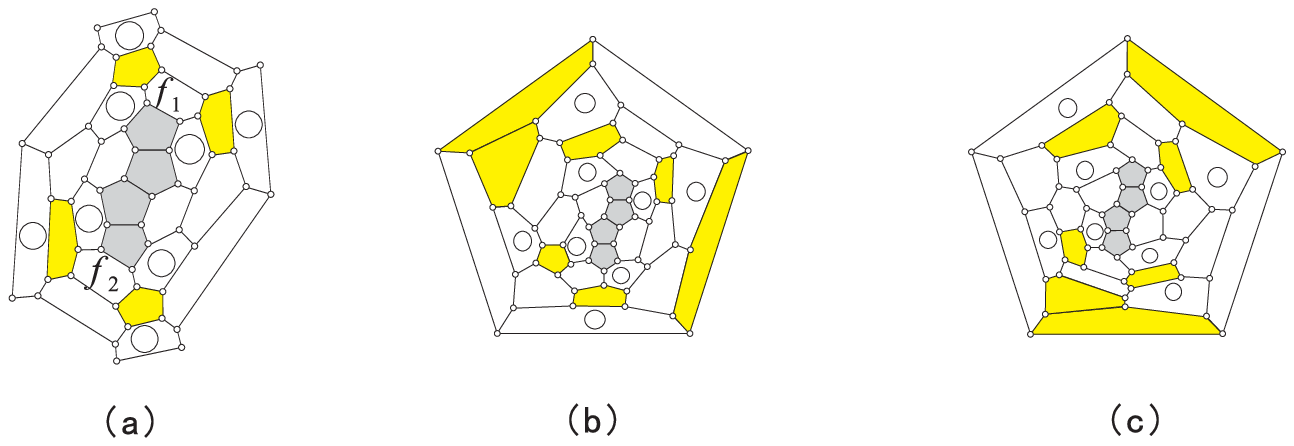}}\\
{Figure \ref{fig4-21}: Illustration for the proof of Case 1 and
extremal fullerene graphs $F^{13}_{60}$ and $F^{14}_{60}$.}
\end{center}
\end{figure}

So suppose $e_1,e_2\in E(B_{1})$ by the symmetry of $e_1$ and $e_2$.
Let $B'_1$ and $B''_1$ be two different subgraphs isomorphic to
$B_1$ such that $e_1\in E(B_1')$ and $e_2\in E(B_1'')$. By
Proposition \ref{pro4-2}, $C[B_1']\cap C[B_1'']=\emptyset$. Hence
$C[B_1']\cup C[B_1'']\cup C[B_2]$ induces the graph (a) in Figure
\ref{fig4-21}. So the remaining four pentagons not in $C[B_1']\cup
C[B_1'']\cup C[B_2]$ adjoin at most four hexagons in $\mathcal H$.
By Lemma \ref{lem4-4}, these four pentagons belong to a $B_{2}$. So
we have two extremal fullerene graphs $F^{13}_{60}$ and
$F^{14}_{60}$ (see Figure \ref{fig4-21} (b) and (c)).

\noindent{\em Case 2:} $e_1\in M$ and $e_2\in E(\mathcal H)$ by
symmetry of $e_1$ and $e_2$. By the discussion of Case 1, we may
assume $e_1\in E(B_1)\cap M$.

\begin{figure}[!hbtp]\refstepcounter{figure}\label{fig4-22}
\begin{center}
\scalebox{0.9}{\includegraphics{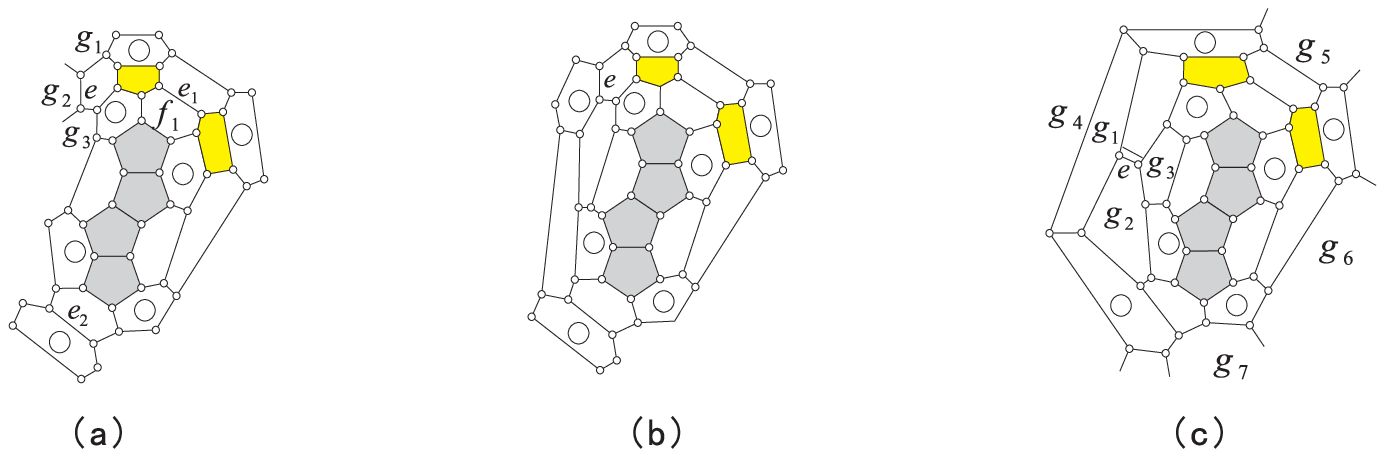}}\\
{Figure \ref{fig4-22}: Illustration for the proof of Case 2.}
\end{center}
\end{figure}

Since every $B_{1}\subset F_{60}$ is maximal, we have the subgraph
of $F_{60}$ as illustrated in Figure \ref{fig4-22} (a). Let $e$ be
an edge on the boundary of the subgraph (a) as shown in Figure
\ref{fig4-22} (a). Then either $e\in E(\mathcal H)$ or $e\in M$. Let
$g_1,g_2,g_3$ be the faces adjoining the subgraph (a) and meeting
$e$. If $e\in E(\mathcal H)$, then $g_2\in \mathcal H$. Hence we
have the graph (b) in Figure \ref{fig4-22}. If the graph (b) is a
subgraph of $F_{60}$, then the remaining six pentagons not in the
graph (b) adjoin at most 5 hexagons in $\mathcal H$, contradicting
Lemma \ref{lem4-4}. So suppose $e\in M$. Then $g_1,g_3$ are
pentagons. Then $g_2$ has to be a hexagon. Hence $e\in E(B_1)\cap
M$. So we have the graph (c) in Figure \ref{fig4-22}. Let
$g_4,g_5,g_6$ and $g_7$ be the faces adjoining the subgraph (c)
along 3-length degree-saturated paths. Note that the graph
consisting of the graph (c) together with $g_4,...,g_7$ has at most
four 2-degree vertices on its boundary. Hence a fullerene graph
$F_n$ containing it satisfies $n\le 58$. So $e_1\in M$ and $e_2\in
E(\mathcal H)$ cannot hold simultaneously.

\noindent{\em Case 3:} $e_1,e_2\in E(\mathcal H)$. By Proposition
\ref{pro4-2}, then $e_1$ and $e_2$ belong to two hexagons in
$\mathcal H$ different from the hexagons in the $C[B_2]$. Let $f_3,
f_4$ be two faces meeting $e_1$ and $e_2$, respectively (see Figure
\ref{fig4-23} (a)).

\begin{figure}[!hbtp]\refstepcounter{figure}\label{fig4-23}
\begin{center}
\scalebox{0.92}{\includegraphics{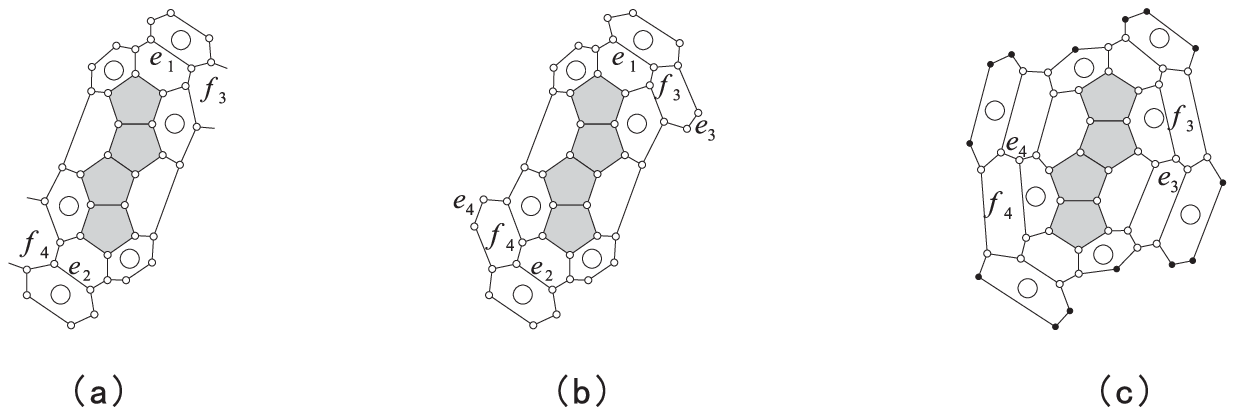}}\\
{Figure \ref{fig4-23}: Illustration for the proof of Case 3.1.}
\end{center}
\end{figure}

\noindent{\em Subcase 3.1:} Both of $f_3$ and $f_4$ are hexagons.
Let $e_3\in E(f_3)$ and $e_4\in E(f_4)$. By Proposition
\ref{pro4-2}, $e_3\ne e_4$ and they are not edges of the graph (a)
(see Figure \ref{fig4-23} (b)).

If $e_3,e_4\in E(\mathcal H)$, then $e_3,e_4$ belong to two distinct
hexagons in $\mathcal H$ and different from the hexagons in the
graph (b). Hence we have the graph (c) in Figure \ref{fig4-23}. The
boundary labeling of the boundary of the graph (c) is 33113113311311
which cannot be separated into the subsequences 13331 (corresponding
to $C_{B_1}-M$) and 1313311 (corresponding to $C_{B_{2}}-M$). So the
graph (c) is not a subgraph of $F_{60}$.

\begin{figure}[!hbtp]\refstepcounter{figure}\label{fig4-24}
\begin{center}
\scalebox{0.9}{\includegraphics{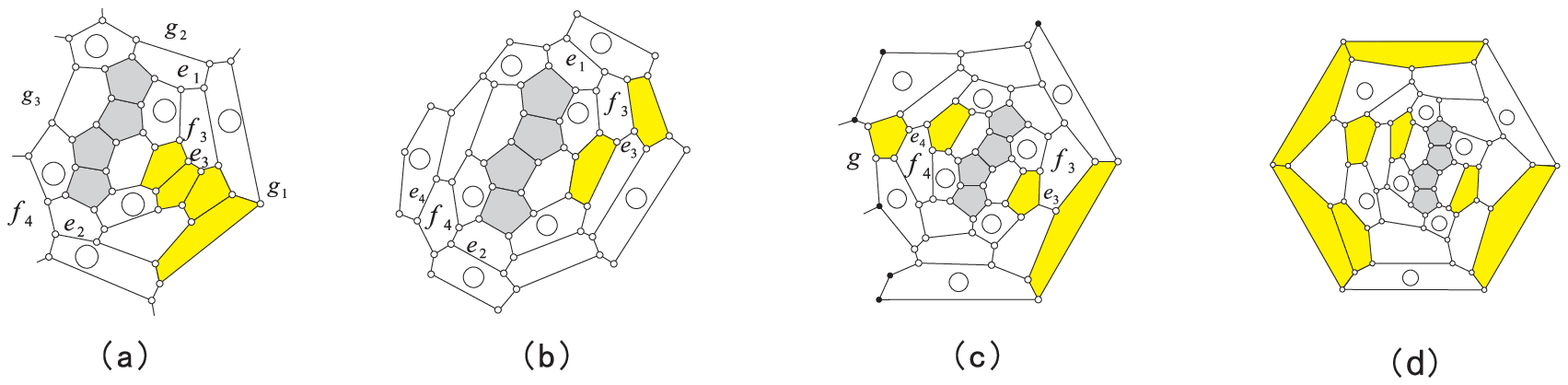}}\\
{Figure \ref{fig4-24}: Illustration for the proof of Case 3.1 and
the extremal fullerene graph $F^{15}_{60}$.}
\end{center}
\end{figure}

So at least one of $e_3$ and $e_4$ belongs to $M$, say $e_3$. If
$e_3\in E(B_{2})$, then we have the graph (a) in Figure
\ref{fig4-24}. Let $g_1,g_2$ and $g_3$ be the faces adjoining it as
shown in Figure \ref{fig4-24} (a). Then the graph consisting of the
graph (a) together with $g_1,g_2,g_3$ and $f_4$ has at most four
2-degree vertices on its boundary. So a fullerene graph $F_n$
containing it satisfies that $n\le 52$ by Lemma \ref{lem2-2}. So
suppose $e_3\in E(B_1)$. If $e_4\in E(\mathcal H)$, then we have the
graph (b) in Figure \ref{fig4-24}. If the graph (b) is a subgraph of
$F_{60}$, then the remaining six pentagons not in the graph (b)
adjoin at most 5 hexagons in $\mathcal H$, contradicting Lemma
\ref{lem4-4}. Therefore, by the symmetry of $e_3$ and $e_4$, we may
assume that $e_4\in E(B_{1})\cap M$. So we have a graph (c) in
Figure \ref{fig4-24}. Since every subgraph of isomorphic to $B_1$ or
$B_{2}$ in $F_{60}$ are maximal, by Lemma \ref{lem4-4}, there is a
unique extremal fullerene graph $F^{15}_{60}$ as shown in Figure
\ref{fig4-24} (d).

\begin{figure}[!hbtp]\refstepcounter{figure}\label{fig4-25}
\begin{center}
\scalebox{0.9}{\includegraphics{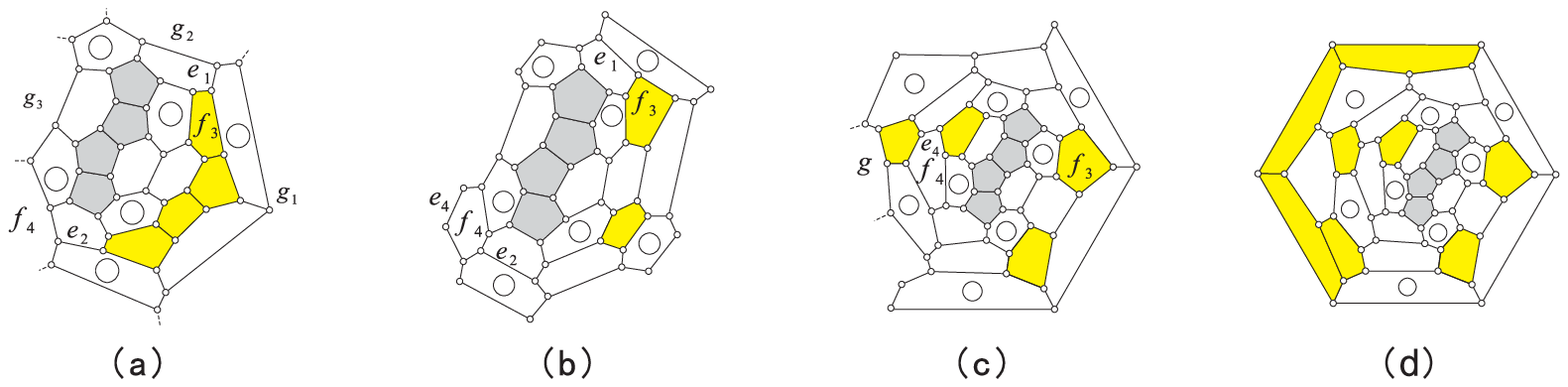}}\\
{Figure \ref{fig4-25}: Illustration for the proof of Subcase 3.2 and
the extremal fullerene graph $F^{16}_{60}$.}
\end{center}
\end{figure}

\noindent{\em Subcase 3.2:} One of $f_3$ and $f_4$ is a pentagon,
say $f_3$. Let $e_4\in E(f_4)$ as that in Subcase 3.1. If $f_3$ is a
pentagon of a $B_{2}$, then we have the graph (a) in Figure
\ref{fig4-25}. A fullerene graph containing the graph (a) has at
most 52 vertices. So suppose $f_3$ is a pentagon of a $B_1$. If
$e_4\in E(\mathcal H)$, then we have a graph (b) in Figure
\ref{fig4-25}. As that $F_{60}$ does not contain the graph (b) in
Figure \ref{fig4-24}, the graph (b) in Figure \ref{fig4-25} is also
not a subgraph of $F_{60}$. Hence $e_4\in M\cap E(B_1)$. Therefore
we have the graph (c) in Figure \ref{fig4-25}. So there is a unique
extremal fullerene graph $F^{16}_{60}$ containing the graph (c)
since every $B_1\subset F_{60}$ is maximal (see Figure \ref{fig4-25}
(d)).

\begin{figure}[!hbtp]\refstepcounter{figure}\label{fig4-26}
\begin{center}
\scalebox{0.95}{\includegraphics{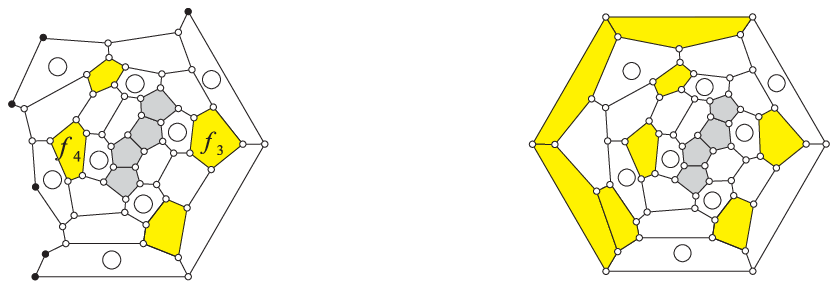}}\\
{Figure \ref{fig4-26}: Illustration for the proof of Subcase 3.3 and
the extremal fullerene graph $F^{17}_{60}$.}
\end{center}
\end{figure}

\noindent{\em Subcase 3.3:} Both $f_3$ and $f_4$ are pentagonal.
According to Subcase 3.2, $f_3$ and $f_4$ belong to subgraphs
isomorphic to $B_1$. Hence, we have a graph as shown in Figure
\ref{fig4-26} (left). Clearly, there are two distinct extremal
fullerene graphs $F_{60}$ containing it: $F^{13}_{60}$ (the graph
(b) in Figure \ref{fig4-21}) and $F^{17}_{60}$ (the right graph in
Figure \ref{fig4-26}).

Combining Cases 1, 2 and 3, there are exact six extremal fullerene
graphs $F_{60}$ such that any $B_1\subset F_{60}$ and any
$B_{2}\subset F_{60}$ are maximal.\end{proof}

\begin{figure}[!hbtp]\refstepcounter{figure}\label{fig4-27}
\begin{center}
\scalebox{1.1}{\includegraphics{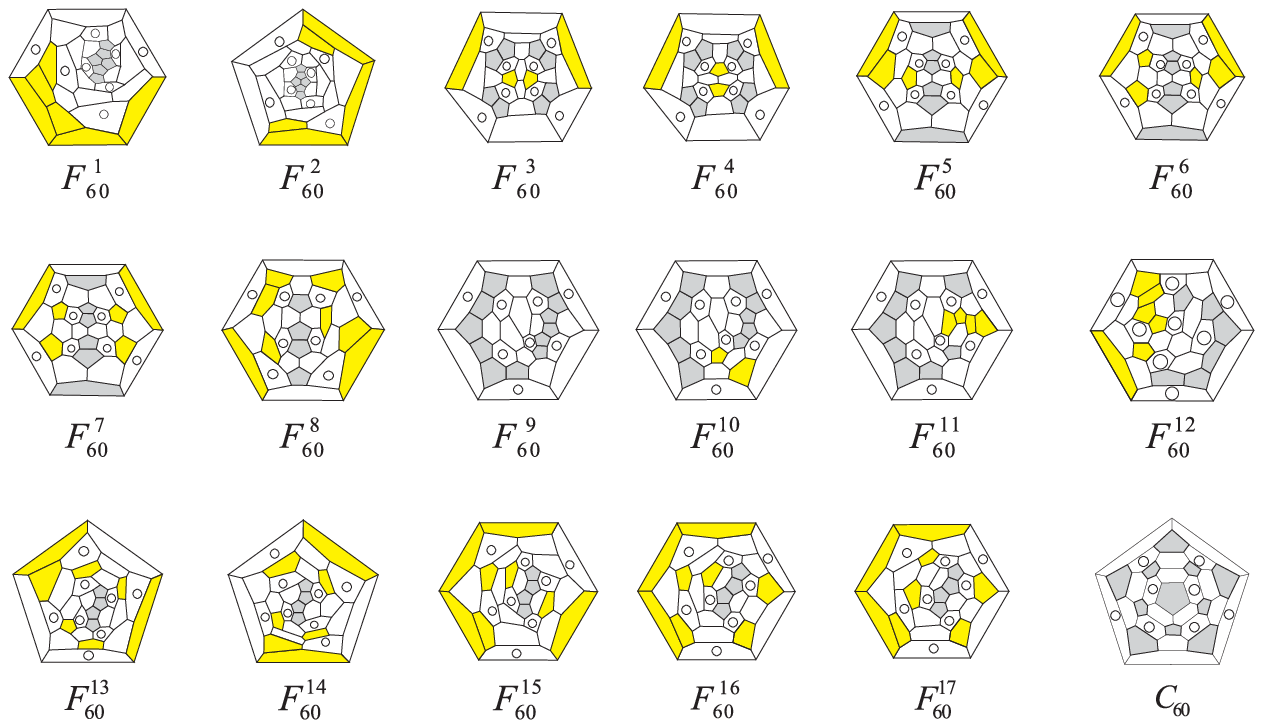}}\\
{Figure \ref{fig4-27}: All extremal fullerene graphs with 60
vertices.}
\end{center}
\end{figure}

Summarizing Lemmas \ref{lem4-6}, \ref{lem4-7}, \ref{lem4-8} and
\ref{lem4-9}, we have the following theorem:

\begin{thm}
There are exactly 18 distinct extremal fullerene graphs with 60
vertices: $\text{C}_{60}$ and $F^i_{60}$ for $i=1,2,...,17$ as shown
in Figure \ref{fig4-27}.
\end{thm}

\section*{Acknowledgements}
The authors are grateful to the referees for their careful reading
and many valuable suggestions.
%%References----------------------------------------------------

\end{document}